\documentclass[a4paper,11pt]{article}
\usepackage{latexsym,amsfonts,amsmath}
\usepackage{graphics}

\def\inter{\mathop{\cap}}
\def\NN{\mathbb{N}}

\def\RR{\mathbb{R}}

\newcommand{\widebar}[1]{\overline{#1}}

\def\union{\mathop{\cup}}
\def\inter{\mathop{\cap}}
\def\liminf{\mathop{\underline{\lim}}}
\def\limsup{\mathop{\overline{\lim}}}
\def\ds{\displaystyle}
\def\SS{\scriptscriptstyle}
\def\nl{\mbox{} \newline }

\usepackage{natbib}

\begin{document}
\newtheorem{theorem}{Theorem}[section]
\newtheorem{proposition}[theorem]{Proposition}
\newtheorem{lemma}[theorem]{Lemma}
\newtheorem{corollary}[theorem]{Corollary}
\newtheorem{definition}[theorem]{Definition}
\newtheorem{remark}[theorem]{Remark}
\newtheorem{conjecture}[theorem]{Conjecture}
\newtheorem{assumption}[theorem]{Assumption}

\bibliographystyle{plain}

\title{Average Continuous Control of Piecewise Deterministic Markov Processes
}
\author{ \mbox{ }
\\
O.L.V. Costa
\thanks{This author received financial support from CNPq (Brazilian National Research
Council), grant 304866/03-2 and FAPESP (Research Council of the
State of S\~ao Paulo), grant 03/06736-7.}
\\ \small Departamento de Engenharia de Telecomunica\c c\~oes e
Controle \\ \small Escola Polit\'ecnica da Universidade de S\~ao
Paulo \\ \small CEP: 05508 900-S\~ao Paulo, Brazil. \\ \small
e-mail: oswaldo@lac.usp.br
\\ \and
\\ F. Dufour
\thanks{Author to whom correspondence should be sent.}
\\
\small Universite Bordeaux I \\
\small IMB, Institut Math\'ematiques de Bordeaux \\
\small INRIA Bordeaux Sud Ouest, Team: CQFD \\
\small \small 351 cours de la Liberation \\
\small 33405 Talence Cedex, France \\ \small e-mail : dufour@math.u-bordeaux1.fr }

\date{ }
\maketitle

\begin{abstract}
This paper deals with the long run average continuous control problem of piecewise deterministic Markov processes (PDMP's)
taking values in a general Borel space and with compact action space depending on the state variable. The control variable acts
on the jump rate and transition measure of the PDMP, and the running and boundary costs are assumed to be positive but not
necessarily bounded. Our first main result is to obtain an optimality equation for the long run average cost in terms of a
discrete-time optimality equation related to the embedded Markov chain given by the post-jump location of the PDMP. Our second main
result guarantees the existence of a feedback measurable selector for the discrete-time optimality equation by establishing
a connection between this equation and an integro-differential equation. Our final main result is to obtain some sufficient
conditions for the existence of a solution for a discrete-time optimality inequality and
an ordinary optimal feedback control for the long run average cost using the so-called
vanishing discount approach (see \cite{hernandez96}, page 83).
\end{abstract}
\begin{tabbing}
\small \hspace*{\parindent}  \= {\bf Keywords:}
piecewise-deterministic Markov Processes, continuous-time,
long-run average\\ cost, optimal control,
integro-differential optimality equation, vanishing approach\\
\> {\bf AMS 2000 subject classification:} \= Primary 60J10
\end{tabbing}

\newpage

\section{Introduction}
A general family of non-diffusion stochastic models suitable for formulating many optimization problems in several areas of
operations research, namely piecewise-deterministic Markov processes (PDMP's), was introduced in \cite{davis84}, and
\cite{davis93}. These processes are determined by three local characteristics; the flow $\phi$, the jump rate $\lambda$ and the
transition measure $Q$. Starting from $x$ the motion of the process follows the flow $\phi(x,t)$ until the first jump time
$T_1$ which occurs either spontaneously in a Poisson-like fashion with rate $\lambda$ or when the flow $\phi(x,t)$ hits the boundary
of the state-space. In either case the location of the process at the jump time $T_1$ is selected by the transition measure
$Q(\phi(x,T_1),.)$ and the motion restarts from this new point as before. A suitable choice of the state space and the local
characteristics $\phi$, $\lambda$, and $Q$ provide stochastic models covering a great number of problems of operations research
\cite{davis93}.

\bigskip

Closely related to the PDMP's are the so called Markov decision drift processes, introduced by Hordijk and Van der Duyn Shouten,
see \cite{hordijk83}, \cite{hordijk84}, \cite{hordijk85}. Their approach is to analyze the control problems for this class of
processes via time-discretizations and deterministic grid approximations of the original process, combined with the theory
of Markov Decision Processes. Yushkevich \cite{yushkevich83,yushkevich87} considers similar control problems but
adopts a variational approach which yields to some integro differential Bellman optimality inequations and characterization
of the value function based on the concept of absolute continuity. By considering some stronger continuity assumptions and under the
framework of the PDMP's, Dempster and Ye\cite{dempster92,dempster96} provide some characterization results
for the value function based on a generalized Bellmam equation which involves the Clark generalized gradient. There exist two
types of control for PDMP's: \textit{continuous control} and \textit{impulse control}. This terminology has been introduced by
M.H.A. Davis in \cite[page 134]{davis93} where continuous control is used to describe situations in which the control variable acts
at all times on the process through the characteristics $(\phi,\lambda,Q)$ by influencing the deterministic motion and the
probability of the jumps by opposition to impulse control that intervenes on the process by moving it to a new point of the state
space at some times specified by the controller.

\bigskip

This paper deals with the long run average continuous control problem of PDMP's taking values in a general Borel space. At each
point $x$ of the state space a control variable is chosen from a compact action set $\mathbb{U}(x)$ and is applied on the jump
parameter $\lambda$ and transition measure $Q$. The goal is to minimize the long run average cost, which is composed of a running
cost and a boundary cost (which is added each time the PDMP touches the boundary). Both costs are assumed to be positive but
not necessarily bounded. As far as the authors are aware of, this is the first time that this kind of problem is considered in the
literature. Indeed, results are available for the long run average cost problem but for impulse control see Costa
\cite{costa89}, Gatarek \cite{garatek93} and the book by M.H.A. Davis \cite{davis93} (see the references therein). On the other
hand, the continuous control problem has been studied only for discounted costs by A. Almudevar \cite{almudevar01}, M.H.A. Davis \cite{davis86,davis93}, M.A.H.
Dempster and J.J. Ye \cite{dempster92,dempster96}, Forwick, Sch{\"a}l, and Schmitz \cite{forwick04}, M. Sch\"al
\cite{schal98}, A.A. Yushkevich \cite{yushkevich87,yushkevich89}. The reader may consult the book by M.H.A. Davis \cite{davis93} and
especially the sections at the end of the chapters 4 and 5 for a complete survey on stochastic control problems for PDMP including
theoretical results and applications.

\bigskip

Our approach to study the long run average continuous control problem of PDMP's  is to follow the idea by M.H.A. Davis and
reformulate the optimal control problem of the PDMP as an equivalent discrete-time Markov decision model in which the stages
are the jump times $T_{n}$. The two main reasons for doing that is to use the powerful tools developed in the discrete-time framework
(see for example the references \cite{bertsekas78,dynkin79,hernandez96,hernandez99}) and to avoid working within the continuous-time context and the associated
infinitesimal generator, which in many situations has its domain difficult to be characterized. For a discounted cost case, the
approach adopted by M.H.A. Davis is very natural since the key idea is to re-write the integral cost as a sum of integrals
between two consecutive jump times of the PDMP and, by doing this, naturally obtaining the one step cost function for the
discrete-time Markov decision model. However, this decomposition for the long run average cost is not possible to be done. Our
first main result is to propose another approach for obtaining an optimality equation for the long run average cost.
It is shown that (see Theorem \ref{theo1a}) if there exist a measurable function $h$, a parameter $\rho$ and a measurable selector satisfying
a discrete-time optimality equation related to the embedded Markov chain given by the post-jump location of the PDMP, and also that an extra condition
involving the function $h$ is verified then an optimal control can be obtained from the measurable selector and $\rho$ is
the optimal cost.

\bigskip

Our second main result (see Theorem \ref{theo1abis}) is to remove the hypothesis of the existence of a measurable selector mentioned in the previous theorem
and in fact, to guarantee the existence of a feedback measurable selector (that is, a selector that depends on the present value of the state
variable, see Remark \ref{feedbak}) provided that the function $h$ and parameter $\rho$ satisfy the optimality equation.
This is done by establishing a link (see the proof Theorem \ref{theo3main}) between the discrete-time optimality equation and an integro-differential equation (using the weaker concept of
absolute continuity along the flow of the value function).
The common approach for the existence of a measurable selector is to impose semicontinuity properties of the cost function and to
introduce the class of relaxed controls to get a compactness property for the action space.
By doing this one obtains an existence result but within the class of relaxed controls. However, what is desired is to show the
existence of an optimal control in the class of ordinary controls. Combining the existence result within the class of relaxed
controls with the connection between the integro-differential equation and the discrete-time equation we can show that the
optimal control is non-relaxed and in fact it is an ordinary feedback control.

In general it is a hard task to get the equality in the solution of the discrete-time optimality equation and verify the extra condition.
A common approach to avoid this is to consider an inequality instead of equality for the optimality equation,
and to use an Abelian result to get the reverse inequality (see for instance \cite{hernandez96}).
Our last main result is to obtain some sufficient conditions, based on the value function of the discounted
control problems, that guarantee the existence of a solution for the discrete-time optimality inequality using the so-called vanishing
approach (see \cite{hernandez96}, page 83). Combining this result with the link between the integro-differential equation and the discrete-time equation
we obtain the existence of an ordinary optimal feedback control for the long run average cost (see Theorem \ref{main}).
In order to do that we need first to establish an optimality equation for the discounted control problem.

A closely related paper to ours, but considering the discounted control case, is the paper by Forwick, Sch{\"a}l, and Schmitz
\cite{forwick04}, which also considers unbounded costs and relaxed controls, and obtain sufficient conditions for the existence of
ordinary feedback controls. However, in \cite{forwick04} the authors do not consider the long run average cost case neither the
related limit problem associated to the vanishing approach. Besides, unlike in \cite{forwick04}, we consider here boundary
jumps and the control action space depending on the state variable. Note however that control on the flow is not considered
here, while it was studied in \cite{forwick04}. Finally it is worth mentioning that the authors are studying in a companion
paper the important question of deriving sufficient stability conditions (like those presented in \cite{costa08},
\cite{dufour99}) under which the conditions on the discounted value function used in the vanishing approach are satisfied,
tracing a parallel with the discrete-time case (see, for instance, \cite{guo06,hernandez96}).

The paper is organized in the following way. In section \ref{pre} we introduce some notation, basic assumptions, and the control
problems to be considered. The definition of the ordinary and relaxed control spaces as well as some operators required for
characterizing the optimality equation are presented in section \ref{conrex}. This section presents several technical
measurability results that will be required throughout the paper. The first main result in presented in section \ref{main1}, Theorem
\ref{theo1a}, which obtains an optimality equation for the long run average
cost in terms of a discrete-time optimality equation related to the embedded Markov chain given by the post-jump location of the PDMP, and an additional
condition. In section \ref{aux1} we introduce some continuity assumptions on the parameters in order to get some convergence and
lower semicontinuity results.
In section \ref{main2} we derive sufficient conditions for the existence of an ordinary feedback optimal control by establishing a link between the discrete-time optimality
equation and an integro-differential equation using the concept of absolute continuity of the value function along the flow (see Theorems \ref{theo1abis}, and \ref{theo3main}).
Section \ref{aux2} considers the discounted optimal control problem and derive an optimality equation.
Our final main result is presented in section \ref{main3} with some sufficient conditions for the existence of a solution for the optimality
inequality and an ordinary optimal feedback control for the long run average cost using the so-called vanishing discount
approach (see Theorem \ref{main}). In order to facilitate the reading of the paper several proofs of some
technical results are presented in the appendix.

\section{Notation and assumptions}
\label{pre} In this section we present some standard notation and some basic definitions related to the motion of a PDMP $\{X(t)\}$,
and the control problems we will consider throughout the paper.
For further details and properties the reader is referred to \cite{davis93}. The following notation will be used in this paper:
\begin{itemize}
\item $\RR$ denotes the set of real numbers, $\RR_+$ the set of positive real numbers and $\RR^d$ the $d$-dimensional euclidian space.
\item $\eta$ denotes the Lebesgue measure on $\RR$.
\item For $X$ a metric spaces, we denote $\mathcal{B}(X)$ as the $\sigma$-algebra generated by the open sets of $X$. $\mathcal{M}(X)$ (respectively, $\mathcal{P}(X)$) denotes the set of all finite (respectively probability) measures on
$(X,\mathcal{B}(X))$. 
\item Let $X$ and $Y$ be metric spaces. The set of all Borel measurable (respectively bounded) functions from $X$ into $Y$ is denoted by $\mathbb{M}(X;Y)$ (respectively $\mathbb{B}(X;Y)$). Moreover, for notational simplicity
$\mathbb{M}(X)$ (respectively $\mathbb{B}(X)$, $\mathbb{M}(X)^{+}$, $\mathbb{B}(X)^{+}$) denotes $\mathbb{M}(X;\RR)$ (respectively $\mathbb{B}(X;\RR)$, $\mathbb{M}(X;\RR_{+})$, $\mathbb{B}(X;\RR_{+})$).
$\mathbb{C}(X)$ denotes the set of continuous functions from $X$ to $\RR$.
For $h\in \mathbb{M}(E)$, $h^+$ (respectively $h^{-}$) denotes the positive (respectively, negtive) part of $h$.
\end{itemize}

\noindent To get a better picture of the motion of a PDMP we first present its definition without any control variable. Let $E$ be an
open subset of $\RR^n$, $\partial E$ its boundary, and $\widebar{E}$ its closure. A PDMP is determined by its local
characteristics $(\phi,\lambda,Q)$. The main assumptions and related definitions on these three parameters are presented below:
\nl
$\bullet$ the flow $\phi(x,t)$ is a function $\phi: \: \mathbb{R}^{n}\times \RR_{+} \longrightarrow \mathbb{R}^{n}$ continuous in 
$(x,t)$ and such that
\begin{eqnarray}\label{pflow}
\phi(x,t+s) & = & \phi(\phi(x,t),s).
\end{eqnarray}
For each $x\in E$ the time the flow takes to reach the boundary starting from $x$ is defined as
$$t_{*}(x)\doteq \inf \{t>0:\phi(x,t)\in \partial E \}.$$

\noindent For $x\in E$ such that $t_{*}(x)=\infty$ (that is, the flow starting from $x$ never touches the boundary), we set
$\phi(x,t_{*}(x))=\Delta$, where $\Delta$ is a fixed point in $\partial E$.

\bigskip

\noindent Some results that will be derived along the paper will be written in terms of the properties along the flow $\phi(x,t)$.
In particular we define the following space of functions absolutely continuous along the flow with limit towards the
boundary:
\begin{align*}
\mathbb{M}^{ac}(E) & =  \bigl\{ g\in\mathbb{M}(E);g(\phi(x,t)): [0,t_{*}(x)) \mapsto \RR \text{ is absolutely continuous for each } x\in E\\
& \text{ and whenever } t_{*}(x)< \infty \text{ the limit }\lim_{t\rightarrow t_{*}(x)} g(\phi(x,t)) \text{ exists} \bigr\}.
\end{align*}
For $g\in \mathbb{M}^{ac}(E)$ and $z\in \partial E$ for which there exists $x\in E$ such that $z=\phi(x,t_{*}(x))$ where $t_{*}(x)< \infty$ we define $\ds g(z) = \lim_{t\rightarrow t_{*}(x)} g(\phi(x,t))$
(note that the limit exists by assumption). As shown in Lemma \ref{Iexist}, for $g\in\mathbb{M}^{ac}(E)$ there exists a function $\mathcal{X}g \in \mathbb{M}(E)$ such that for all $x\in E$ and $t\in [0,t_{*}(x))$
$$g(\phi(x,t))-g(x)  =  \int_{0}^{t} \mathcal{X}g(\phi(x,s)) ds .$$
\nl
$\bullet$ the jump rate $\lambda:E\rightarrow \RR_{+}$ which is assumed to be a measurable function satisfying: $(\forall x \in E)$
$(\exists \varepsilon >0)$ such that $\ds \int_0^\varepsilon \lambda(\phi(x,s))ds< \infty$.
\nl $\bullet$ the post-jump location kernel $Q: \widebar{E} \times \mathcal{B}(E) \rightarrow [0,1]$ which is a transition measure satisfying the following property:
$(\forall x\in E)$ $Q(x,E-\{x\})=1$.

\bigskip

\noindent From these characteristics, it can be shown \cite[p.62-66]{davis93} that there exists a filtered probability space
$(\Omega,\mathcal{F},\{ \mathcal{F}_{t} \}, \{ P_{x} \}_{x\in E})$
such that the motion of the process $\{X(t)\}$ starting from a point $x\in E$ may be constructed as follows. Take a random variable $T_1$ such that
\begin{equation*}
P_x(T_1>t) \doteq
\begin{cases}
e^{-\Lambda(x,t)} & \text{for } t<t_{*}(x)\\
0 & \text{for } t\geq t_{*}(x)
\end{cases}
\end{equation*}
and $t\in [0,t_*(x)[$
\begin{equation*}
\Lambda(x,t) \doteq \int_0^t\lambda(\phi(x,s))ds.
\end{equation*}
If $T_1$ generated according to the above probability is equal to infinity, then for $t\in \RR_+$, $X(t)=\phi(x,t)$. Otherwise select independently an $E$-valued random variable (labelled
$X_{1}$) having distribution $Q(\phi(x,T_1),.)$.
The trajectory of $\{X(t)\}$ starting at $x$, for $t\leq T_1$ , is given by
\begin{equation*}
X(t) \doteq
\begin{cases}
\phi(x,t)&\text{for } t<T_1, \\
X_1&\text{for }t=T_1.
\end{cases}
\end{equation*}
Starting from $X(T_1)=X_1$, we now select the next inter-jump time $T_2-T_1$ and post-jump location $X(T_2)=X_2$ in a similar way.
This gives a strong Markov process $\{X(t)\}$ with jump times $\bigl\{T_k\bigr\}_{k \in \NN}$ (where $T_{0}\doteq 0$).

\bigskip

We present next the control problems and some basic assumptions we will consider throughout the paper. We suppose from now on that the local characteristics $\lambda$ and $Q$ depend on a control
action $u\in \mathbb{U}$ where $\mathbb{U}$ is a Borel space, in the following way:
\begin{itemize}
\item[$\bullet$] $\lambda \in \mathbb{M}(\widebar{E}\times\mathbb{U})^{+}$. \item[$\bullet$] $Q$ is a stochastic kernel on $E$ given $\widebar{E}\times \mathbb{U}$.
\end{itemize}

\noindent For each $x\in \widebar{E}$ we define the subsets $\mathbb{U}(x)$ of $\mathbb{U}$ as the set of feasible control actions that can be taken when the state process is in $x\in \widebar{E}$,
that is, the control action that will be applied to $\lambda$ and $Q$ must belong to $\mathbb{U}(x)$.
The following assumptions, based on the standard theory of Markov decision processes (see \cite{hernandez96}), will be made throughout the paper:
\begin{assumption}
\label{Hyp1a} For all $x\in \widebar{E}$, $\mathbb{U}(x)$ is a compact subspace of $\mathbb{U}$.
\end{assumption}
\begin{assumption}
\label{Hyp2a} The set $K=\left\{(x,a): x\in \widebar{E}, a \in \mathbb{U}(x) \right\}$ is a Borel subset of $\widebar{E}\times \mathbb{U}$.
\end{assumption}

\noindent The following assumption will also be required along the paper:
\begin{assumption}
\label{Hyp4a} For all $x\in E$, and $t\in [0,t_{*}(x))$, $\ds
\int_{0}^{t} \sup_{a\in \mathbb{U}(\phi(x,s))}
\lambda(\phi(x,s),a) \: ds < \infty$. \nl If $t_{*}(x)< \infty$
then $\ds \int_{0}^{t_{*}(x)} \sup_{a\in \mathbb{U}(\phi(x,s))}
\lambda(\phi(x,s),a) \: ds < \infty$.
\end{assumption}

\noindent We present next the definition of an admissible control strategy and the associated motion of the controlled process.
A control policy $U$ is a pair of functions $(u,u_{\partial}) \in \mathbb{M}(\NN \times E\times \RR_{+};\mathbb{U}) \times \mathbb{M}(\NN \times E;\mathbb{U})$
satisfying $u(n,x,t)\in \mathbb{U}(\phi(x,t))$, and $u_{\partial}(n,x)\in \mathbb{U}(\phi(x,t_{*}(x)))$ for all $(n,x,t)\in \NN \times E\times \RR_{+}$.
\nl
The class of admissible control strategy will be denoted by $\mathcal{U}$.

\bigskip

\noindent Given a control strategy $U=(u,u_{\partial})$, one describe the motion of the piecewise deterministic process $X(t)$ in the following manner. Define $T_{0}=0$ and $X(0)=x$.
Assume that the process $\{X(t)\}$ is located at $Z_{n}$ at the $n^{th}$ jump time $T_{n}$ then select a random variable $S_{n}$ having distribution
$$F(t)=1-I_{\{t<t_{*}(Z_{n})\}} e^{-\int_{0}^{t}\lambda(\phi(Z_{n},t),u(n,Z_{n},s)) ds}.$$
Define $T_{n+1}=T_{n}+S_{n}$ and for $t\in[T_{n},T_{n+1})$, $X(t)=\phi(Z_{n},t-T_{n})$. \nl Let $Z_{n+1}$ a random variable having distribution $Q(\phi(Z_{n},T_{n+1}),u(n,Z_{n},S_{n}));.)$
if $\phi(Z_{n},T_{n+1})\in E$ or $Q(\phi(Z_{n},T_{n+1}),u_{\partial}(n,Z_{n});.)$ if $\phi(Z_{n},T_{n+1})\in\partial E$. \nl At time $T_{n+1}$, the process $\{X(t)\}$ is defined by $X(T_{n+1})=Z_{n+1}$.

\bigskip

\noindent Now we give a more precise definition of the controlled piecewise deterministic Markov process described above. Consider
the state space $\widehat{E}=E\times E \times \RR_{+} \times \NN$.
For a control policy $U=(u,u_{\partial})$ let us introduce the following parameters for $\hat{x}=(x,z,s,n)\in \widehat{E}$:
\nl
$\bullet$ the flow $\widehat{\phi}(\hat{x},t) = (\phi(x,t), z , s+t , n)$. \nl $\bullet$ the jump rate $\widehat{\lambda}^{U}(\hat{x})=\lambda(x,u(n,z,s))$.
\nl
$\bullet$ the transition measure
\begin{eqnarray*}
\widehat{Q}^{U}(\hat{x},A\times B\times \{0\}\times \{n+1\}) =
\begin{cases}
Q(x,u(n,z,s)); A\inter B) & \text{ if } x\in E,\\
Q(x,u_{\partial}(n,z);A\inter B) & \text{ if } x\in \partial E,
\end{cases}
\end{eqnarray*}
for $A$ and $B$ in $\mathcal{B}(E)$.

\bigskip

\noindent From \cite[section 25]{davis93}, it can be shown that for any control strategy $U=(u,u_{\partial})\in
\mathcal{U}$ there exists a filtered probability space $(\Omega,\mathcal{F},\{ \mathcal{F}_{t} \}, \{ P^{U}_{\hat{x}}
\}_{\hat{x}\in \widehat{E}})$ such that the piecewise deterministic Markov process $\{\widehat{X}^{U}(t)\}$ with local
characteristics $(\widehat{\phi},\widehat{\lambda}^{U},\widehat{Q}^{U})$ may be constructed as follows. For notational simplicity the probability
$P^{U}_{\hat{x}_{0}}$ will be denoted by $P^{U}_{(x,k)}$ for $\hat{x}_{0}=(x,x,0,k)\in \widehat{E}$. Take a random variable
$T_1$ such that
\begin{equation*}
P^{U}_{(x,k)}(T_1>t) \doteq
\begin{cases}
e^{-\Lambda^{U}(x,k,t)} & \text{for } t<t_{*}(x)\\
0 & \text{for } t\geq t_{*}(x)
\end{cases}
\end{equation*}
where for $x\in E$ and $t\in [0,t_*(x)[$
\begin{equation*}
\Lambda^{U}(x,k,t) \doteq \int_0^t\lambda(\phi(x,s),u(k,x,s))ds.
\end{equation*}
If $T_1$ is equal to infinity, then for $t\in \RR_+$,
$\widehat{X}^{U}(t)= \bigl(\phi(x,t),x,t,0\bigr)$. Otherwise
select independently an $\widehat{E}$-valued random variable
(labelled $\widehat{X}^{U}_{1}$) having distribution
\begin{eqnarray*}
P^{U}_{(x,k)}(\widehat{X}^{U}_{1} \in A \times B \times \{0\}
\times \{k+1\} |\sigma\{T_1\}) =
\begin{cases}
Q(\phi(x,T_{1}),u(k,x,T_{1})); A\inter B) \text{ if } \phi(x,T_{1})\in E, \\
Q(\phi(x,T_{1}),u_{\partial}(k,x); A\inter B) \text{ if } \phi(x,T_{1}) \in \partial E.
\end{cases}
\end{eqnarray*}
The trajectory of $\{\widehat{X}^{U}(t)\}$ starting from
$(x,x,0,k)$, for $t\leq T_1$ , is given by
\begin{equation*}
\widehat{X}^{U}(t) \doteq
\begin{cases}
\bigl(\phi(x,t),x,t,0\bigr) &\text{for } t<T_1, \\
\widehat{X}^{U}_1 &\text{for } t=T_1.
\end{cases}
\end{equation*}
Starting from $\widehat{X}^{U}(T_1)=\widehat{X}^{U}_1$, we now
select the next inter-jump time $T_2-T_1$ and post-jump location
$\widehat{X}^{U}(T_2)=\widehat{X}^{U}_2$ in a similar way.

\bigskip

\noindent Let us define the components of the PDMP
$\{\widehat{X}^{U}(t)\}$ by
$\widehat{X}^{U}(t)=\bigl(X(t),Z(t),\tau(t),N(t)\bigr)$. From the
previous construction of the PDMP $\{\widehat{X}^{U}(t)\}$, it is
easy to see that $X(t)$ corresponds to the trajectory of the
system, $Z(t)$ is the value of $X(t)$ at the last jump time before
$t$, $\tau(t)$ is time elapsed between the last jump and time $t$,
and $N(t)$ is the number of jumps of the process  $\{X(t)\}$ at
time $t$.

\bigskip

\noindent As in Davis \cite{davis93}, we consider the following
assumption to avoid any accumulation point of the jump times:
\begin{assumption}
\label{Hypjump} For any $x\in E$, $U=(u,u_{\partial})\in
\mathcal{U}$, and $t\geq 0$, $\ds E^{U}_{(x,0)}\Biggl[
\sum_{i=1}^{\infty} I_{\{T_{i}\leq t\}} \Biggr] < \infty$
\end{assumption}

\noindent The costs of our control problem will contain two terms, a running cost $f$ and a boundary cost $r$, satisfying the following properties:
\begin{assumption}
\label{Hyp6a} $f\in \mathbb{M}(\widebar{E}\times\mathbb{U})^{+}$.
\end{assumption}
\begin{assumption}
\label{Hyp7a} $r\in \mathbb{M}(\partial E\times\mathbb{U})^{+}$.
\end{assumption}

\noindent The long-run average cost we want to minimize over
$\mathcal{U}$ is given by:
\begin{align}
\mathcal{A}(U,x) = \limsup_{t\rightarrow +\infty} \frac{1}{t} E^{U}_{(x,0)} \Biggl[ \int_{0}^{t} f\bigl(X(s), & u(N(s),Z(s),\tau(s)) \bigr) ds \nonumber \\
& + \int_{0}^{t} r\bigl(X(s-),u_{\partial}(N(s-),Z(s-)) \bigr)
dp^{*}(s) \Biggr], \label{avcost}
\end{align}
where $\ds p^{*}(t) = \sum_{i=1}^{\infty} I_{\{T_{i}\leq t\}}
I_{\{ X(T_{i}-) \in \partial E\}}$ counts the number of times the
process hits the boundary up to time $t$, and we set
\begin{align}
\mathcal{J}_{\mathcal{A}}(x) =
\inf_{U\in\mathcal{U}}\mathcal{A}(U,x). \label{avcostvf}
\end{align}

\bigskip

\noindent For the $\alpha$ discounted case, with $\alpha > 0$, the
cost we want to minimize is given by:
\begin{align}
\mathcal{D}^{\alpha}(U,x) = E^{U}_{(x,0)} \Biggl[ \int_{0}^{\infty} e^{-\alpha s} f\bigl(X(s), & u(N(s),Z(s),\tau(s)) \bigr) ds \nonumber \\
& + \int_{0}^{\infty} e^{-\alpha s}
r\bigl(X(s-),u_{\partial}(N(s-),Z(s-)) \bigr) dp^{*}(s) \Biggr]
\label{discost}
\end{align}
and we set
\begin{align}
\mathcal{J}_{\mathcal{D}}^{\alpha}(x) =
\inf_{U\in\mathcal{U}}\mathcal{D}^{\alpha}(U,x). \label{discostvf}
\end{align}

\bigskip

\noindent We also consider a truncated version of problem
(\ref{discostvf}) defined, for each $m=0,1,\ldots$, as
\begin{align}
\mathcal{D}^{\alpha}_{m}(U,x) = E^{U}_{(x,0)} \Biggl[ \int_{0}^{T_m} e^{-\alpha s} f\bigl(X(s), & u(N(s),Z(s),\tau(s)) \bigr) ds \nonumber \\
& + \int_{0}^{T_m} e^{-\alpha s}
r\bigl(X(s-),u_{\partial}(N(s-),Z(s-)) \bigr) dp^{*}(s) \Biggr]
\label{discostm}
\end{align}

\bigskip

\noindent We need the following assumption, to avoid infinite
costs for the discounted case.
\begin{assumption}
\label{Hypdis} For all $\alpha>0$ and all $x\in E$,
$\mathcal{J}_{\mathcal{D}}^{\alpha}(x)<\infty$.
\end{assumption}

It is clear that for all $x\in E$, $\ds 0\leq
\inf_{U\in\mathcal{U}}\mathcal{D}^{\alpha}_{m}(U,x)\leq
\mathcal{J}_{\mathcal{D}}^{\alpha}(x)<\infty$.

\begin{remark}
There is no loss of generality in assuming that $\mathbb{U}$ is compact.
Indeed, if $\mathbb{U}$ is a Borel space, it follows from Proposition 7 in \cite{yushkevich97} that $\mathbb{U}$ can be considered as a measurable subset of a compact space $\mathbb{U}'$
where $\mathbb{U}(x)$ is compact in $\mathbb{U}'$ for $x\in \widebar{E}$.
Moreover, by recalling that $\mathbb{U}(x)$ represents the set of feasible controls in the state $x\in \widebar{E}$,
the definition of the control problem will not be affected if the functions $\lambda$, $Q$, $f$, and $r$ are extended to $\mathbb{U}'$.
Therefore, from now on one will assume that $\mathbb{U}$ is compact. This result will be needed in sub-section \ref{conrex1}
\end{remark}

\section{Discrete-time ordinary and relaxed controls}
\label{conrex} The class of open loop policies denoted by $\mathcal{U}$ has been introduced in the previous section as the set of admissible control strategies for a PDMP.
As mentioned in the introduction we will study in section \ref{main1} how the original continuous-time control problem can be associated to an optimality equation of a discrete-time problem related to the embedded Markov
chain given by the post-jump location of the PDMP.
In this section we first present the definitions of the discrete-time ordinary and relaxed control sets used in the formulation of the optimality equation of the discrete-time Markov control problem as well as
the characterization of some topological properties of these sets.
In particular, by using a result of the theory of multifunctions (see the book by Castaing and Valadier \cite{castaing77}), it is shown that the set of relaxed controls is compact.
In the sequel we present some important operators associated to the optimality equation of the discrete-time problem as well as some measurability properties.

\subsection{Relaxed and ordinary control}
\label{conrex1}
We present in this sub-section the set of discrete-time relaxed controls and the subset of ordinary controls. Consider the Banach
spaces $L^{1}(\RR_{+};\mathbb{C}(\mathbb{U}))$ and $L^{\infty}(\RR_{+};\mathcal{M}(\mathbb{U}))$ where
$\mathbb{C}(\mathbb{U})$ is equipped with the topology of uniform convergence and $\mathcal{M}(\mathbb{U})$ is equipped with the
weak$^{*}$ topology $\sigma(\mathcal{M}(\mathbb{U}),\mathbb{C}(\mathbb{U}))$. Let $\mathcal{V}^{r}$ (respectively $\mathcal{V}^{r}(x)$ for $x\in E$)
be the set of all $\eta$-measurable functions $\mu$ defined on $\RR_{+}$ with value in $\mathcal{P}(\mathbb{U})$ such that
$\mu(t,\mathbb{U})=1$ $\eta$-a.e. (respectively $\mu(t,\mathbb{U}(\phi(x,t)))=1$ $\eta$-a.e.). From Theorem V-2 in
\cite{castaing77}, it follows that $\mathcal{V}^{r}$ (respectively $\mathcal{V}^{r}(x)$ for $x\in E$) are compact sets with respect
to the weak$^{*}$ topology $\sigma(L^{\infty}(\RR_{+};\mathcal{M}(\mathbb{U})),L^{1}(\RR_{+};\mathbb{C}(\mathbb{U})))$.
Moreover, from Bishop's Theorem (see Theorem I.3.11 in \cite{warga72}), there is a metric such for all $x\in E$,
$\mathcal{V}^{r}(x)$ is a compact set of the Borel set $\mathcal{V}^{r}$. Note that a sequence $\bigl(\mu_{n}\bigr)_{n\in
\NN}$ in $\mathcal{V}^{r}(x)$ converges to $\mu$ if and only if \begin{eqnarray}
\lim_{n\rightarrow \infty} \int_{\RR_{+}} \int_{\mathbb{U}(\phi(x,t))}  g(t,u) \mu_{n}(t,du) dt =
\int_{\RR_{+}} \int_{\mathbb{U}(\phi(x,t))}  g(t,u) \mu(t,du) dt, \label{convergence}
\end{eqnarray}
for all $g\in L^{1}(\RR_{+};\mathbb{C}(\mathbb{U}))$.

\bigskip

\noindent Therefore, the set of relaxed controls are defined as follows. For $x\in E$,
\begin{align*}
\mathbb{V}^{r}(x) & = \mathcal{V}^r(x) \times \mathcal{P}\bigl(\mathbb{U}(\phi(x,t_{*}(x)))\bigr), \\
\mathbb{V}^{r} & = \mathcal{V}^r \times
\mathcal{P}\bigl(\mathbb{U}\bigr).
\end{align*}
The set of ordinary controls, denoted by $\mathbb{V}$ (respectively $\mathbb{V}(x)$ for $x\in E$), is defined as above
except that it is composed of deterministic functions instead of probability measures. More specifically we have
\begin{eqnarray*}
\mathcal{V}(x) & = & \Bigl\{ \nu \in \mathbb{M}(\RR_{+}, \mathbb{U}) : (\forall t\in \RR_{+}), \nu(t) \in \mathbb{U}(\phi(x,t)) \Bigr\}, \\
\mathbb{V}(x) & = & \mathcal{V}(x) \times \mathbb{U}(\phi(x,t_{*}(x))), \\
\mathbb{V} & = & \mathbb{M}(\RR_{+}, \mathbb{U}) \times
\mathbb{U}.
\end{eqnarray*}

\noindent Consequently, the set of ordinary controls is a subset of the set of relaxed controls $\mathbb{V}^{r}$ (respectively
$\mathbb{V}^{r}(x)$ for $x\in E$) by identifying any control action $u\in \mathbb{U}$ with the Dirac measure concentrated on
$u$. Thus we can write that $\mathbb{V}\subset \mathbb{V}^{r}$ (respectively $\mathbb{V}(x)\subset \mathbb{V}^{r}(x)$ for $x\in
E$) and from now on we will consider that $\mathbb{V}$ (respectively $\mathbb{V}(x)$ for $x\in E$) will be endowed with
the topology generated by $\mathbb{V}^{r}$.

\bigskip

\noindent The necessity to introduce the class of relaxed control $\mathbb{V}^{r}$ is justified by the fact that in general there
does not exist a topology for which $\mathbb{V}$ and $\mathbb{V}(x)$ are compact sets. However from the previous
construction, it follows that $\mathbb{V}^{r}$ and $\mathbb{V}^{r}(x)$ are compact sets.

\bigskip

\noindent We present next a definition that will be useful in the next sections.
\begin{definition}
\label{def1*} For any $x\in E$, $t\in [0,t_{*}(x))$, and
$\Theta=\bigl( \mu, \mu_{\partial} \bigr)\in \mathbb{V}^{r}(x)$,
define
\begin{eqnarray}
\bigl[\Theta\bigr]_{t}=\bigl( \mu(.+t), \mu_{\partial} \bigr).
\label{DefUps-rt}
\end{eqnarray}
Clearly,  $\bigl[\Theta\bigr]_{t} \in \mathbb{V}^{r}(\phi(x,t))$.
\end{definition}

\bigskip

\noindent As in \cite{hernandez96}, page 14, we need that the set of feasible state/relaxed-control pairs is a measurable subset of $\mathcal{B}(E)\times \mathcal{B}(\mathbb{V}^{r})$, that is, we need the following assumption.
\begin{assumption}
\label{Mesurability}
$$\mathcal{K} \doteq \bigl\{ (x,\Theta) : \Theta \in \mathbb{V}^{r}(x), x\in E \bigr\} \in \mathcal{B}(E)\times \mathcal{B}(\mathbb{V}^{r})$$
\end{assumption}

\bigskip

\noindent We present a sufficient condition, based on the continuity of the sets $\mathbb{U}(x)$, to ensure that assumption
\ref{Mesurability} holds. The proof is presented in the appendix.
\begin{proposition}\label{mesura1}
Assumption \ref{Mesurability} is satisfied if for all convergent sequence $\{x_{n}\}_{n\in \NN}$ in $E$
$$\mathbb{U}(x)=\inter_{n\in \NN} \union_{m\geq n} \mathbb{U}(x_{m})=  \union_{n\in \NN} \inter_{m\geq n} \mathbb{U}(x_{m}),$$
where  $\ds \lim_{n\rightarrow \infty} x_{n} = x.$
\end{proposition}

\subsection{Discrete-time operators and measurability properties}
In this sub-section we present some important operators associated to the optimality equation of the discrete-time problem as well as
some measurability properties.

We consider the following notation for $x\in \widebar{E}$, $y\in E$ and
$\mu\in\mathcal{P}\bigl(\mathbb{U}\bigr)$, $h\in \mathbb{M}(E)^{+}$, and $w\in \mathbb{M}(\widebar{E}\times \mathbb{U})^{+}$:
\begin{align*}
& Qh(x,\mu) =  \int_{\mathbb{U}} \int_{E} h(z) Q(x,u;dz) \mu (du), \qquad \lambda Qh(x,\mu) =  \int_{\mathbb{U}} \lambda(x,u) \int_{E} h(z) Q(x,u;dz)  \mu (du),\\
&\Lambda^{\mu}(y,t)  =  \int_{0}^{t} \int_{\mathbb{U}} \lambda(\phi(y,s),u) \mu (du) ds ,\,\,\,\,\,\qquad w(x,\mu)  =  \int_{\mathbb{U}} w(x,u)  \mu (du).&
\end{align*}

\noindent The following operators will be associated to the optimality equations of the discrete-time problems that will be
presented in the next sections. For $\Theta=\bigl(\mu,\mu_{\partial}\bigr)\in \mathbb{V}^{r}$, $h\in \mathbb{M}(E)$, $\alpha \in \RR$, let us introduce the following stochastic kernel defined on $E\times \mathcal{B}(E)$ by
\begin{eqnarray}
\label{DefGr}
G_{\alpha}(x,\Theta;A) & \doteq & \int_0^{t_{*}(x)}e^{-\alpha s - \Lambda^{\mu}(x,s)}\lambda QI_{A}(\phi(x,s),\mu(s)) ds\nonumber \\
& &+ e^{-\alpha t_{*}(x) -\Lambda^{\mu}(x,t_{*}(x))} Q(\phi(x,t_{*}(x)),\mu_{\partial};A),
\end{eqnarray}
for all $(x,A)\in E\times \mathcal{B}(E)$, and so for $h\in \mathbb{M}(E)^{+}$, we define $G_{\alpha}h(x,\Theta)\doteq\ds \int_{E} h(y) G_{\alpha}(x,\Theta;dy)$.
For $x\in E$, $\Theta=\bigl(\mu,\mu_{\partial}\bigr)\in \mathbb{V}^{r}$, $v\in \mathbb{M}(E\times \mathbb{U})^{+}$, $w\in \mathbb{M}(\partial E\times \mathbb{U})$,
$\alpha \in \RR$, introduce
\begin{eqnarray}
\label{DefLr}
L_{\alpha}v(x,\Theta) & \doteq  & \int_0^{t_{*}(x)}e^{-\alpha s-\Lambda^{\mu}(x,s)} v(\phi(x,s),\mu(s)) ds, \\
\label{DefHr} H_{\alpha}w(x,\Theta) & \doteq & e^{-\alpha t_{*}(x)-\Lambda^{\mu}(x,t_{*}(x))} w(\phi(x,t_{*}(x)),\mu_{\partial}).
\end{eqnarray}
For $h\in \mathbb{M}(E)$ (respectively, $v\in \mathbb{M}(E\times \mathbb{U})$),  $G_{\alpha}h(x,\Theta)=G_{\alpha}h^{+}(x,\Theta)-G_{\alpha}h^{-}(x,\Theta)$
(respectively, $L_{\alpha}v(x,\Theta)=L_{\alpha}v^{+}(x,\Theta)-L_{\alpha}v^{-}(x,\Theta)$) provided the difference has a meaning.
It will be useful in the sequel to define the function $\mathcal{L}_{\alpha}(x,\Theta)$  as follows:
$\mathcal{L}_{\alpha}(x,\Theta)  \doteq  L_{\alpha}I_{E\times \mathbb{U}}(x,\Theta)$. In particular for $\alpha =0$ we write for simplicity $G_0=G$, $L_0=L$, $H_0=H$, $\mathcal{L}_0=\mathcal{L}$.

\bigskip

\noindent The next proposition presents some important measurability properties of the operators $G_{\alpha}$,
$L_{\alpha}$, and $H_{\alpha}$ (defined in equations (\ref{DefGr}), (\ref{DefLr}), and (\ref{DefHr})) and its proof can
be found in the appendix.

\begin{proposition}
\label{mes3} Let $\alpha\in \RR_{+}$, $g\in \mathbb{M}(E)$ be bounded from below, $w_{1}\in \mathbb{M}(E\times \mathbb{U})$ be
bounded from below, and $w_{2}\in \mathbb{M}(\partial E\times \mathbb{U})$. Then the mappings $G_{\alpha}g(x,\Theta)$,
$L_{\alpha}w_{1}(x,\Theta)$, and $H_{\alpha}w_{2}(x,\Theta)$ defined on $E\times \mathbb{V}^r$ with values in $\RR$ are
$\mathcal{B}(E\times \mathbb{V}^r)$-measurable.
\end{proposition}

\bigskip

\noindent We present now the definitions of the one-stage optimization operators.
\begin{definition}
For $\alpha\in \RR_{+}$, $\rho\in \RR_{+}$, and $g\in \mathbb{M}(E)$
bounded from below, define
\begin{itemize}
\item[i)] the (ordinary) one-stage optimization operator by
\begin{eqnarray}
\mathcal{T}_{\alpha}(\rho,g)(x) & = & \inf_{\Upsilon\in \mathbb{V}(x)} \Bigl\{-\rho \mathcal{L}_{\alpha}(x,\Upsilon) +
L_{\alpha}f(x,\Upsilon)+H_{\alpha}r(x,\Upsilon)+G_{\alpha} g(x,\Upsilon)  \Bigr\}.
\end{eqnarray}
\item[ii)] the relaxed one-stage optimization operator by
\begin{eqnarray}
\mathcal{R}_{\alpha}(\rho,g)(x) & = & \inf_{\Theta\in \mathbb{V}^{r}(x)} \Bigl\{-\rho \mathcal{L}_{\alpha}(x,\Theta) +
L_{\alpha}f(x,\Theta)+H_{\alpha}r(x,\Theta)+G_{\alpha} g(x,\Theta) \Bigr\}.
\end{eqnarray}
\end{itemize}
\end{definition}

\noindent In particular for $\alpha =0$ we write for simplicity $\mathcal{T}_{0}=\mathcal{T}$, and $\mathcal{R}_{0}=\mathcal{R}$.

\bigskip

\noindent Let us introduce the following sets of measurable selectors associated to $\bigl(\mathbb{U}(x)\bigr)_{x\in E}$
(respectively $\bigl(\mathbb{V}(x)\bigr)_{x\in E}$, $\bigl(\mathbb{V}^{r}(x)\bigr)_{x\in E}$):
\begin{eqnarray*}
\mathcal{S}_{\mathbb{U}} & = & \Bigl\{  u \in \mathbb{M}(\widebar{E}, \mathbb{U}) : (\forall x\in \widebar{E}), u(x) \in \mathbb{U}(x)\Bigr\}, \\
\mathcal{S}_{\mathbb{V}} & = & \Bigl\{ (\nu,\nu_{\partial})\in \mathbb{M}(E, \mathbb{V}) : (\forall x\in E), \bigl( \nu(x), \: \nu_{\partial}(x) \bigr)\in \mathbb{V}(x)\Bigr\}, \\
\mathcal{S}_{\mathbb{V}^{r}} & = &\Bigl\{ (\mu,\mu_{\partial})\in
\mathbb{M}(E, \mathbb{V}^{r}) : (\forall x\in E), \bigl( \mu(x),
\: \mu_{\partial}(x) \bigr)\in \mathbb{V}^{r}(x)\Bigr\}.
\end{eqnarray*}

\begin{remark}
The set $\mathcal{S}_{\mathbb{U}}$ characterizes the control law $u(x)$ that only depends on the value of the state variable $x$.
On the other hand $(\nu,\nu_{\partial})\in \mathcal{S}_{\mathbb{V}}$ characterizes an ordinary control law for the control problem associated to
the one-stage optimization operator.
Indeed, starting from $x$, it defines the control law for all $t\in[0,t_*(x))$ through the function $\nu(x,t)$ and at $t=t_*(x)$ (if $t_*(x)<\infty$) through $\nu_{\partial}$.
Finally $(\mu,\mu_{\partial})\in \mathcal{S}_{\mathbb{V}^{r}}$
characterizes a relaxed control law for the control problem associated to the relaxed one-stage optimization operator.
Since starting from $x$, it defines a probability over the feasible control actions for all $t\in[0,t_*(x))$ through the probability measure $\mu(x,t)$ and at $t=t_*(x)$ (if $t_*(x)<\infty$) through the probability
measure $\mu_{\partial}$.
\end{remark}

\bigskip

\noindent For $\alpha\in \RR_{+}$, $\rho\in \RR_{+}$, and $g\in \mathbb{M}(E)$ bounded from below, the one-stage optimization problem associated to the operator $\mathcal{T}_{\alpha}(\rho,g)$, respectively
$\mathcal{R}_{\alpha}(\rho,g)$, consists of finding a measurable selector $\Upsilon\in \mathcal{S}_{\mathbb{V}}$, respectively $\Theta\in \mathcal{S}_{\mathbb{V}^{r}}$ such that for all $x\in E$
\begin{eqnarray*}
\mathcal{T}_{\alpha}(\rho,g)(x) & = & -\rho \mathcal{L}_{\alpha}(x,\Upsilon(x)) + L_{\alpha}f(x,\Upsilon(x))+H_{\alpha}r(x,\Upsilon(x))+G_{\alpha} g(x,\Upsilon(x)),
\end{eqnarray*}
respectively
\begin{eqnarray*}
\mathcal{R}_{\alpha}(\rho,g)(x) & = & -\rho \mathcal{L}_{\alpha}(x,\Theta(x)) + L_{\alpha}f(x,\Theta(x))+H_{\alpha}r(x,\Theta(x))+G_{\alpha} g(x,\Theta(x)).
\end{eqnarray*}

\bigskip

\noindent Finally we conclude this section by showing that there exist two natural mappings from $\mathcal{S}_{\mathbb{U}}$ to $\mathcal{S}_{\mathbb{V}}$
and from $\mathcal{S}_{\mathbb{U}}$ to $\mathcal{U}$
\begin{definition}
\label{mapu}
For $u \in\mathcal{S}_{\mathbb{U}}$, define the mapping $u_{\phi}$ $:$ $x$ $\rightarrow$
$\bigl(u(\phi(x,.)),u(\phi(x,t_{*}(x)))\bigr)$ of the space $E$ into $\mathbb{V}$
\end{definition}

\begin{proposition}
\label{lem3} If $u \in\mathcal{S}_{\mathbb{U}}$ then $u_{\phi} \in \mathcal{S}_{\mathbb{V}}$.
\end{proposition}
\textbf{Proof:} From Lemma A.3 in \cite{forwick04} and item $(i)$ of Lemma 3 in \cite{yushkevich80} it follows that the mapping $x$
$\rightarrow$ $u(\phi(x,.))$ of the space $E$ into $\mathbb{M}(\RR_{+}, \mathbb{U})$ is measurable. Moreover, for all
$(x,t)\in E\times \RR_{+}$, $u(\phi(x,t)) \in \mathbb{U}(\phi(x,t))$. Therefore, $u_{\phi}$ belongs to
$\mathcal{S}_{\mathbb{V}}$. \hfill $\Box$

\begin{definition}
\label{mapU}
For $u \in\mathcal{S}_{\mathbb{U}}$, define the mapping $U_{\phi}$ $:$ $(n,x,t)$ $\rightarrow$
$\bigl(u(\phi(x,t)),u(\phi(x,t_{*}(x)))\bigr)$ of the space $\NN\times E\times \RR_{+}$ into $\mathbb{U}\times \mathbb{U}$.
\end{definition}

\begin{proposition}
\label{lem4} If $u \in\mathcal{S}_{\mathbb{U}}$ then $U_{\phi} \in \mathcal{U}$.
\end{proposition}
\textbf{Proof:} This is a straightforward consequence of the measurability properties of $u$ and $\phi$ and the fact that $u(x)\in \mathbb{U}(x)$.
\hfill $\Box$

\begin{remark}
\label{feedbak} The measurable selectors of the kind $u_{\phi}$ as in Definition \ref{mapu} are called feedback measurable selectors in the class $\mathcal{S}_{\mathbb{V}}\subset \mathcal{S}_{\mathbb{V}^{r}}$
and the control strategies of the kind $U_{\phi}$ as in definition \ref{mapU} are called feedback control strategies in the class $\mathcal{U}$.
\end{remark}

\section{Discrete-time optimality equation for the average control}
\label{main1} In this section we obtain an optimality equation for
the long run average cost problem defined in equation
(\ref{avcostvf})
in terms of a discrete-time optimality equation related to the embedded Markov chain given by the post-jump location of the PDMP, and an additional 
condition on a limit over the solution of the optimality equation
divided by the time $t$.
Notice that in this section we will be assuming that there is an optimal solution
for the one-stage optimization problem, and therefore we will be dealing with the ordinary action sets (that is, $\mathbb{V}(x)$).

\bigskip

\begin{theorem}
\label{theo1a} Suppose that there exists a pair $(\rho,h) \in \RR_{+}\times \mathbb{M}(E)$ with $h$ bounded from below satisfying
the following discrete-time optimality equation
\begin{eqnarray}
\mathcal{T}(\rho,h)(x)  & = & h(x), \label{OptEq0}
\end{eqnarray}
and for all $U\in \mathcal{U}$,
\begin{eqnarray}
\limsup_{t\rightarrow +\infty} \frac{1}{t} \limsup_{m\rightarrow +\infty } E^{U}_{(x,0)} \Bigl[  h\bigl(X(t\wedge T_{m})
\bigr)\Bigr] & = & 0. \label{oeavcost2}
\end{eqnarray}
Moreover, assume that there exists a solution to the one-stage optimization problem associated to $\mathcal{T}(\rho,h)$, that is, the existence of an optimal measurable selector
$\hat{\Gamma}=(\hat{\gamma}, \hat{\gamma}_{\partial})$ in $\mathcal{S}_{\mathbb{V}}$ such that for all $x\in E$
\begin{eqnarray}
\mathcal{T}(\rho,h)(x) & = & -\rho \mathcal{L}(x,\hat{\Gamma}(x)) + Lf(x,\hat{\Gamma}(x))+Hr(x,\hat{\Gamma}(x))+Gh(x,\hat{\Gamma}(x)),
\label{oeavcost}
\end{eqnarray}
Define the control strategy $\widehat{U}$ by $(\hat{u},\hat{u}_{\partial})$ with $\hat{u}(n,x,t)=\hat{\gamma}(x,t)$, $\hat{u}_{\partial}(n,x)=\hat{\gamma}_{\partial}(x)$ for $(n,x,t)\in \NN \times E \times \RR_{+}$.
Then $\widehat{U}$ belongs to $\mathcal{U}$ and it is optimal. Moreover,
\begin{eqnarray*}
\rho  = \mathcal{J}_{\mathcal{A}}(x) = \mathcal{A}(\widehat{U},x).
\end{eqnarray*}
\end{theorem}
The proof of this theorem is presented at the end of this section.

\bigskip

\noindent In Theorem \ref{theo1a} notice that equation (\ref{OptEq0}) can be seen as the optimality equation of a
discrete-time problem related to the embedded Markov chain given by the post-jump location of the PDMP with transition kernel $G$, and equation (\ref{oeavcost2}) as an additional technical condition.

\bigskip

\noindent Notice that in general, one cannot guarantee the existence of an optimal measurable selector for the optimality
equation (\ref{OptEq0}) without compactness conditions.
This problem of existence will be considered in the next section.

\bigskip

\noindent In order to prove the previous theorem we first need to present several intermediate results that will also be used in the remainder of the paper.
Notice that these results will be written in terms of an extra parameter $\alpha\geq 0$ that will be useful for the discounted control problem analyzed in
sections \ref{aux2} and \ref{main3}.
The proofs of these intermediate results can be found in the appendix.

\bigskip

\noindent The next proposition presents some important properties of the one-stage optimality equation (see equations (\ref{OpDr1a}), (\ref{Defw*})).
It is shown that the solution of the one-stage optimality equation has a special time representation (see equation (\ref{eqOpDr1b})).
As a consequence it follows that it is absolutely continuous along trajectories with limit on the boundary (that is, it belongs to $\mathbb{M}^{ac}(E)$).

\begin{proposition}
\label{prop3a} Let $\rho\in \RR_{+}$ and $h\in\mathbb{M}(E)$ be bounded from below. For $\alpha\geq 0$ and  $x\in E$ define
\begin{equation}
w(x)=\mathcal{T}_{\alpha}(\rho,h)(x). \label{OpDr1a}
\end{equation}
Assume that $w\in \mathbb{M}(E)$ and there exists $\hat{\Gamma}\in \mathcal{S}_{\mathbb{V}}$ such that
\begin{eqnarray}
w(x)=-\rho \mathcal{L}_{\alpha}(x,\hat{\Gamma}(x)) + L_{\alpha}f(x,\hat{\Gamma}(x))
+H_{\alpha}r(x,\hat{\Gamma}(x))+G_{\alpha}h(x,\hat{\Gamma}(x)).
\label{Defw*}
\end{eqnarray}
Then $w\in\mathbb{M}^{ac}(E)$ and for all $x\in E$, and $t\in [0,t_{*}(x))$,
\begin{align}
w(x)= & \int_0^t e^{-\alpha s-\Lambda^{\hat{\gamma}(x)}(x,s)}\biggl[ -\rho + f(\phi(x,s),\hat{\gamma}(x,s)) + \lambda
Qh(\phi(x,s),\hat{\gamma}(x,s))\biggr] ds
\nonumber \\
& + e^{-\alpha t-\Lambda^{\hat{\gamma}(x)}(x,t)}w(\phi(x,t)) \label{eqOpDr1b}\\
= & \inf_{ \nu \in \mathcal{V}(x)} \bigg\{ \int_0^t e^{-\alpha s-\Lambda^{\nu}(x,s)} \biggl[ -\rho + f(\phi(x,s),\nu(s)) + \lambda Qh(\phi(x,s),\nu(s))\biggr] ds \nonumber \\
& + e^{-\alpha t-\Lambda^{\nu}(x,t)}w(\phi(x,t))\bigg\},
\label{eqOpDr1a}
\end{align}
where $\hat{\Gamma}(x)=\bigl( \hat{\gamma}(x),
\hat{\gamma}_{\partial}(x) \bigr)$.
\end{proposition}

\bigskip

\noindent
The next two propositions deal with two inequalities of opposite directions for the one-stage optimality equation.
Roughly speaking these two results show that if $h$ is a solution for a one-stage optimality inequality (see (\ref{OpD3a}) or (\ref{OpD3b})) then this inequality is
preserved, in one case for any control strategy and in the other case for a specific control strategy, along the jump time
iterations for a cost conveniently defined, see equations (\ref{iqJma}) or (\ref{iqJmb}).

\begin{proposition}
\label{prop2a}
Let $h\in\mathbb{M}(E)$ be bounded from below. For $\rho \in \RR_{+}$ and $\alpha\in \RR_{+}$,
assume that $\mathcal{T}_{\alpha}(\rho,h) \in \mathbb{M}(E)$
and there exists $\hat{\Gamma}=(\hat{\gamma}, \hat{\gamma}_{\partial}) \in \mathcal{S}_{\mathbb{V}}$ such that for all $x\in E$
\begin{align}
h(x) & \leq \mathcal{T}_{\alpha}(\rho,h)(x) =  -\rho \mathcal{L}_{\alpha}(x,\hat{\Gamma}(x)) +
L_{\alpha}f(x,\hat{\Gamma}(x))+H_{\alpha}r(x,\hat{\Gamma}(x))+G_{\alpha}h(x,\hat{\Gamma}(x)).
\label{OpD3a}
\end{align}
For $U\in \mathcal{U}$, $m\in \NN$, and $(t,x,k)\in \RR_{+}\times E\times \NN$, define
\begin{align*}
J^{U}_{m}(t,x,k) =  & E^{U}_{(x,k)} \biggl[ \int_{0}^{t\wedge T_{m}} e^{-\alpha s}\Bigl[ f\bigl(X(s),u(N(s),Z(s),\tau(s)) \bigr) -\rho \Bigr] ds \nonumber \\
& + \int_{0}^{t\wedge T_{m}}e^{-\alpha s}r\bigl(X(s-),u_{\partial}(N(s),X(s-)) \bigr) dp^{*}(s) +
e^{-\alpha (t\wedge T_{m})} \mathcal{T}_{\alpha}(\rho,h)\bigl(X(t\wedge T_{m}) \bigr)\biggr].
\end{align*}
Then for all $m\in \NN$, and $(t,x,k)\in \RR_{+}\times E\times
\NN$,
\begin{align}
J^{U}_{m}(t,x,k)  \geq h(x). \label{iqJma}
\end{align}
\end{proposition}

\bigskip

\noindent The next proposition considers the reverse inequality.
\begin{proposition}
\label{prop2b} Let $h\in\mathbb{M}(E)$ be bounded from below. For $\rho \in \RR_{+}$ and $\alpha\in \RR_{+}$,
assume that $\mathcal{T}_{\alpha}(\rho,h) \in \mathbb{M}(E)$
and there exists $\hat{\Gamma}=(\hat{\gamma}, \hat{\gamma}_{\partial}) \in \mathcal{S}_{\mathbb{V}}$ such that for all $x\in E$
\begin{align}
h(x) & \geq \mathcal{T}_{\alpha}(\rho,h)(x) =  -\rho \mathcal{L}_{\alpha}(x,\hat{\Gamma}(x)) + L_{\alpha}f(x,\hat{\Gamma}(x))+H_{\alpha}r(x,\hat{\Gamma}(x))+G_{\alpha}h(x,\hat{\Gamma}(x)).
\label{OpD3b}
\end{align}
Then for $\widehat{U}$ defined by $(\hat{u},\hat{u}_{\partial})$ with $\hat{u}(n,x,t)=\hat{\gamma}(x,t)$, and
$\hat{u}_{\partial}(n,x)=\hat{\gamma}_{\partial}(x)$ for $(n,x,t)\in \NN \times E \times \RR_{+}$, we have that
$\widehat{U}$ belongs to $\mathcal{U}$. Moreover, defining
\begin{align*}
J^{\widehat{U}}_{m}(t,x,k) =  & E^{\widehat{U}}_{(x,k)} \biggl[ \int_{0}^{t\wedge T_{m}} e^{-\alpha s}\Bigl[ f\bigl(X(s),\hat{u}(N(s),Z(s),\tau(s)) \bigr) -\rho \Bigr] ds \nonumber \\
& + \int_{0}^{t\wedge T_{m}}e^{-\alpha s}
r\bigl(X(s-),\hat{u}_{\partial}(N(s),X(s-)) \bigr) dp^{*}(s)
+ e^{-\alpha (t\wedge T_{m})} \mathcal{T}_{\alpha}(\rho,h)\bigl(X(t\wedge T_{m}) \bigr)\biggr],
\end{align*}
we have, for all $m\in \NN$, $(t,x,k)\in \RR_{+}\times E\times \NN$, that
\begin{align}
J^{\widehat{U}}_{m}(t,x,k)  \leq h(x). \label{iqJmb}
\end{align}
\end{proposition}

\bigskip
\noindent Combining the previous propositions with $\alpha=0$ we
get the proof of Theorem \ref{theo1a}.
\bigskip

\noindent {\textbf{Proof of Theorem \ref{theo1a}}:} From
Proposition \ref{prop2a} it follows that
\begin{align*}
E^{U}_{(x,0)} \biggl[ & \int_{0}^{t\wedge T_{m}} f\bigl(X(s),u(N(s),Z(s),\tau(s)) \bigr) ds  + \int_{0}^{t\wedge T_{m}} r\bigl(X(s-),u_{\partial}(N(s),X(s-)) \bigr) dp^{*}(s)\biggr] \nonumber \\
& + E^{U}_{(x,0)} \Bigl[ h\bigl(X(t\wedge T_{m}) \bigr)\Bigr] \geq  E^{U}_{(x,0)} \Bigl[ \rho [t\wedge T_{m}] \Bigr] +h(x).
\end{align*}
Consequently we have, from assumption \ref{Hypjump} (which implies
that $T_{m} \rightarrow \infty$ $P^{U}$ a.s.), that
\begin{align*}
\limsup_{t\rightarrow +\infty} \frac{1}{t} E^{U}_{(x,0)} \biggl[ & \int_{0}^{t} f\bigl(X(s),u(N(s),Z(s),\tau(s)) \bigr) ds  + \int_{0}^{t} r\bigl(X(s-),u_{\partial}(N(s),X(s-)) \bigr) dp^{*}(s)\biggr] \nonumber \\
& +\limsup_{t\rightarrow +\infty} \frac{1}{t} \limsup_{m\rightarrow +\infty } E^{U}_{(x,0)} \Bigl[ h\bigl(X(t\wedge T_{m}) \bigr)\Bigr] \geq  \rho,
\end{align*}
showing that $\rho \leq \mathcal{J}_{\mathcal{A}}(x)$, by using equation (\ref{oeavcost2}).

\bigskip

\noindent From Proposition \ref{prop2b}, the control strategy $\widehat{U}$ defined by $(\hat{u},\hat{u}_{\partial})$ with
$\hat{u}(n,x,t)=\hat{\gamma}(x,t)$, $\hat{u}_{\partial}(n,x)=\hat{\gamma}_{\partial}(x)$ for $(n,x,t)\in \NN \times E \times \RR_{+}$ belongs to $\mathcal{U}$.
Combining equations (\ref{oeavcost2}) and (\ref{iqJmb}) we obtain that $\rho \geq \mathcal{A}(\widehat{U},x)$ completing the proof.
\hfill $\Box$

\section{Convergence and lower semicontinuity properties}
\label{aux1} In the previous section we assumed the existence of an ordinary optimal measurable selector for the one-stage optimization problem associated to $\mathcal{T}(\rho,h)$
(see equation (\ref{oeavcost})), for $(\rho,h)$ satisfying the optimality equation $\mathcal{T}(\rho,h)(x) = h(x)$.
In the next sections we will suppress this hypothesis.
In order to do that we need to consider relaxed controls, so that we can take advantage of the compactness property of the sets $\mathbb{V}^{r}(x)$
and $\mathbb{V}^{r}$ as presented in section \ref{conrex}.
Note however that we also need the cost function to be lower semicontinuous.
Thus in this section we present the assumptions and results that will guarantee some convergence and
lower semicontinuity properties of the operators $G_{\alpha}$, $L_{\alpha}$, and $H_{\alpha}$
that appear in the one-stage optimization operators with respect to the topology defined in equation (\ref{convergence}).
Combining the compactness of the sets $\mathbb{V}^{r}(x)$ with the lower semicontinuity of the operators $G_{\alpha}$, $L_{\alpha}$, and
$H_{\alpha}$ we can use the measurable selector theorem as presented in Proposition D.5 of \cite{hernandez96} to get the
existence of a relaxed optimal control and measurability of the one-stage optimization equation.
Moreover in parallel we get some important convergence properties that will be applied with the vanishing approach in section \ref{main3}.

\bigskip

\noindent From now on we will consider the following assumptions.

\begin{assumption}
\label{Hyp3a} For each $x\in E$, $\lambda(x,.):\mathbb{U}(x)\mapsto \RR_{+}$ is continuous.
\end{assumption}
\begin{assumption}
\label{Hyp6bis} There exists a sequence of measurable functions
$(f_j)_{j\in \NN}$ in $\mathbb{M}(\widebar{E}\times \mathbb{U})^{+}$ such that for all $y\in \widebar{E}$, $f_j(y,.)\uparrow f(y,.)$ as $j\rightarrow \infty$ and $f_j(y,.) \in \mathbb{C}(\mathbb{U}(y))$.
\end{assumption}
\begin{assumption}
\label{Hyp7bis} There exists a sequence of measurable functions
$(r_j)_{j\in \NN}$ in $\mathbb{M}(\partial E\times \mathbb{U})^{+}$ such that for all $z\in \partial E$,
$r_j(z,.)\uparrow r(z,.)$ as $j\rightarrow \infty$ and $r_j(z,.) \in \mathbb{C}(\mathbb{U}(z))$.
\end{assumption}
\begin{assumption}
\label{Hyp5a} For all $x \in \widebar{E}$ and $g\in \mathbb{B}(E)$, $Qg(x,.):\mathbb{U}(x)\mapsto \RR$ is continuous.
\end{assumption}
\begin{assumption}
\label{Hyp8a} There exists $\xi \in\mathbb{M}(E)^{+}$, such that
\begin{enumerate}
\item [a)] $\lambda(y,a) \geq \xi(y)$ for all $y\in E$ and $a\in\mathbb{U}(y)$,
\item [b)] $\ds \int_0^{t_{*}(x)} e^{-\int_0^t \xi(\phi(x,s))ds}dt \leq K_{\xi}< \infty$ for all $x\in E$,
\item [c)] $\ds \int_0^{t_{*}(x)} e^{-\int_0^t \xi(\phi(x,s))ds} \sup_{a\in \mathbb{U}(\phi(x,t))}f(\phi(x,t),a)dt < \infty$.
\end{enumerate}
\end{assumption}

\begin{remark} \label{vac} A consequence of Assumption \ref{Hyp8a} b) is
that $\ds \lim_{t\rightarrow +\infty} e^{-\alpha t-\int_0^{t}\xi(\phi(x,s))ds} =0$, for any $\alpha \in \RR_{+}$ and $x\in E$ with $t_{*}(x)=+\infty$.
Therefore, for any $x\in E$ with $t_{*}(x)=+\infty$, $A\in \mathcal{B}(E)$, $\alpha\in \RR_{+}$, $\Theta=(\mu,\mu_{\partial})\in \mathbb{V}^{r}(x)$,
$w\in \mathbb{M}(\partial E\times \mathbb{U})$, $\ds G_{\alpha}(x,\Theta;A) = \int_0^{t_{*}(x)}e^{-\alpha s - \Lambda^{\mu}(x,s)}\lambda QI_{A}(\phi(x,s),\mu(s)) ds$, and
$H_{\alpha}w(x,\Theta) =0$.
\end{remark}

\noindent The next proposition presents convergence results of the operators $G_{\alpha}$, $L_{\alpha}$, and $H_{\alpha}$ with
respect to the topology defined in equation (\ref{convergence}).
Note that the convergence is taken not only with respect to a sequence of controls but also with respect to some functions and the parameter $\alpha$.
This is justified by the fact we will need this convergence for the vanishing approach in section \ref{main3}.
The proof of the proposition is in the appendix.

\begin{proposition}
\label{lemlsc1} Consider $\alpha \in \RR_{+}$ and a non increasing sequence of positive numbers $\{\alpha_k\}$, $\alpha_k\downarrow \alpha$, a sequence of functions $h_{\alpha_k}\in \mathbb{M}(E)$
uniformly bounded from below by a positive constant $K_h$ (that is, $h_{\alpha_k}(y) \geq -K_h$ for all $y\in E$).
Set $\ds h = \liminf_{k\rightarrow \infty} h_{\alpha_k}$.
For $x\in E$, consider $\Theta_n=\bigl(\mu_n,\mu_{\partial,n}\bigr)\in \mathbb{V}^{r}(x)$ and $\Theta=\bigl(\mu,\mu_{\partial}\bigr)\in \mathbb{V}^{r}(x)$ such that $\Theta_n \rightarrow \Theta$.
We have the following results:
\begin{enumerate}
\item [a)] $\ds \lim_{n\rightarrow \infty}\mathcal{L}_{\alpha_n}(x,\Theta_n)=\mathcal{L}(x,\Theta)$.
\item [b)] $\ds \liminf_{n\rightarrow \infty}L_{\alpha_n}f(x,\Theta_n)\geq Lf(x,\Theta)$.
\item [c)] $\ds \liminf_{n\rightarrow \infty}H_{\alpha_n}r(x,\Theta_n)\geq Hr(x,\Theta)$. \item [d)] $\ds \liminf_{n\rightarrow \infty} G_{\alpha_n}h_{\alpha_n}(x,\Theta_n) \geq Gh(x,\Theta)$.
\end{enumerate}
\end{proposition}

\noindent The lower semicontinuity properties mentioned at the beginning of this section follow easily from this proposition as stated in the next corollary.

\begin{corollary}
\label{corlsc1} Consider $h \in \mathcal{M}(E)$  bounded from below. We have the following results:
\begin{enumerate}
\item [a)] $\mathcal{L}_\alpha(x,\Theta)$ is continuous on $\mathbb{V}^{r}(x)$.
\item [b)] $L_{\alpha} f(x,\Theta)$ is lower semicontinuous on $\mathbb{V}^{r}(x)$.
\item [c)] $H_{\alpha} r(x,\Theta)$ is lower semicontinuous on $\mathbb{V}^{r}(x)$.
\item [d)] $G_{\alpha} h(x,\Theta)$ is lower semicontinuous on $\mathbb{V}^{r}(x)$.
\end{enumerate}
\end{corollary}
\textbf{Proof:} By taking $\alpha_k=\alpha\geq 0$, $h_{\alpha_k}=h$ in Proposition \ref{lemlsc1} the results follow.
\hfill $\Box$

\section{Existence of an ordinary optimal feedback control}
\label{main2}
The main result of this section is the Theorem \ref{theo1abis} that strengthens Theorem \ref{theo1a} of the previous section by
only assuming that the discrete-time optimality equation $\mathcal{T}(\rho,h) =  h$ has a solution in order to ensure the existence of an optimal control strategy
for the long run average control problem. Moreover, it is shown that this optimal control strategy is in the feedback class and can be characterized as in item $D3)$ of the Definition \ref{defmes}.
These results are obtained by establishing a connection (see the proof Theorem \ref{theo3main}) between the discrete-time optimality equation and an integro-differential equation (using the weaker concept of
absolute continuity along the flow of the value function).
The basic idea is to use the set of relaxed controls $\mathbb{V}^{r}(x)$. The advantage of considering
$\mathbb{V}^{r}(x)$ is that it is compact so that, together with the assumptions we have made in section \ref{aux1},
we can apply a measurable selector theorem to guarantee the existence of an optimal measurable selector (see Proposition
\ref{prop3b}). The price to pay is that this measurable selector belongs to the space of relaxed controls.
However  we can show that in fact there exists a non-relaxed feedback selector for the discrete-time optimality equation $\mathcal{T}(\rho,h) =  h$ by establishing
a connection between the discrete-time optimality equation and the integro-differential equation (see the proof of Theorem \ref{theo3main}).

\bigskip

\begin{definition}
\label{defmes}
Consider $w\in\mathbb{M}(E)$ and $h\in\mathbb{M}(E)$ bounded from below.
\begin{enumerate}
\item [D1)] Denote by $\widehat{u}(w,h)\in \mathcal{S}_{\mathbb{U}}$ the measurable selector satisfying
\begin{align*}
\inf_{a\in  \mathbb{U}(x)} \{f(x,a) -\lambda(x, & a)  \Bigl[ w(x)-Q h(x,a) \Bigr]\} \\
& =  f(x,\widehat{u}(w,h)(x))-\lambda(x,\widehat{u}(w,h)(x)) \Bigl[w(x)-Qh(x,\widehat{u}(w,h)(x)) \Bigr],
\end{align*}
\begin{eqnarray*}
\inf_{a\in \mathbb{U}(z)}\{r(z,a)+Qh(z,a)\} & = & r(z,\widehat{u}(w,h)(z))+Qh(z,\widehat{u}(w,h)(z)).
\end{eqnarray*}
\item [D2)] $\widehat{u}_{\phi}(w,h)\in \mathcal{S}_{\mathbb{V}}$ is the measurable selector derived from $\widehat{u}(w,h)$ through
the definition \ref{mapu}.
\item [D3)] $\widehat{U}_{\phi}(w,h)\in \mathcal{U}$ is the control strategy derived from $\widehat{u}(w,h)$ through
the definition \ref{mapU}.
\end{enumerate}
\end{definition}
The existence of $\widehat{u}(w,h)$ follows from assumptions \ref{Hyp3a}-\ref{Hyp5a} and Theorem 3.3.5 in \cite{hernandez96} and
the fact that $\widehat{u}_{\phi}(w,h)\in \mathcal{S}_{\mathbb{V}}$, and $\widehat{U}_{\phi}(w,h)\in \mathcal{U}$
comes from Propositions \ref{lem3} and \ref{lem4}.

\bigskip

\begin{theorem}
\label{theo1abis} Suppose that there exists a pair $(\rho,h) \in \RR\times \mathbb{M}(E)$ with $h$ bounded from below satisfying
the following discrete-time optimality equation
\begin{eqnarray*}
\mathcal{T}(\rho,h)(x) & = & h(x) ,
\end{eqnarray*}
and for all $U\in \mathcal{U}$,
\begin{eqnarray*}
\limsup_{t\rightarrow +\infty} \frac{1}{t} \limsup_{m\rightarrow +\infty } E^{U}_{(x,0)} \Bigl[  h\bigl(X(t\wedge T_{m})
\bigr)\Bigr] & = & 0.
\end{eqnarray*}
Then $h\in\mathbb{M}^{ac}(E)$, the feedback optimal control strategy $\widehat{U}_{\phi}(h,h)$ (see item D3) of Definition \ref{defmes})
is optimal, and
\begin{eqnarray*}
\rho  = \mathcal{J}_{\mathcal{A}}(x) = \mathcal{A}(\widehat{U}(h,h),x).
\end{eqnarray*}
\end{theorem}
\textbf{Proof:} The proof of this result is straightforward by combining Theorem \ref{theo1a} of the previous section and Theorem \ref{theo3main} presented below.
 \hfill $\Box$

\bigskip

\noindent The proof of the next proposition is presented in the appendix.
It shows the existence of an optimal relaxed measurable selector for the relaxed one-stage optimization operator $\mathcal{R}_{\alpha}(\rho,h)(x)$ and that $\mathcal{R}_{\alpha}(\rho,h) \in\mathbb{M}^{ac}(E)$.
\begin{proposition}
\label{prop3b}
Let $\alpha\geq 0$, $\rho \in \RR_{+}$ and $h\in\mathbb{M}(E)$ be bounded from below. For $x\in E$ define
$w(x)=\mathcal{R}_{\alpha}(\rho,h)(x)$. Assume that for all $x\in E$, $w(x)\in \RR$.
Then there exists $\hat{\Theta}\in \mathcal{S}_{\mathbb{V}^{r}}$ such that
\begin{align}
w(x) & = -\rho \mathcal{L}_{\alpha}(x,\hat{\Theta}(x)) + L_{\alpha}f(x,\hat{\Theta}(x))+H_{\alpha}r(x,\hat{\Theta}(x))
+G_{\alpha} h(x,\hat{\Theta}(x)).
\label{O2prop3b}
\end{align}
Moreover, $w\in\mathbb{M}^{ac}(E)$, and satisfies for all $x\in E$ and $t\in [0,t_{*}(x))$
\begin{align}
w(x)  & = \inf_{ \mu \in \mathcal{V}^{r}(x)} \bigg\{ \int_0^t e^{-\alpha s-\Lambda^{\mu}(x,s)}\biggl[ -\rho + f(\phi(x,s),\mu(s))
+ \lambda Qh(\phi(x,s),\mu(s))\biggr] ds \nonumber \\
& \phantom{=} + e^{-\alpha t-\Lambda^{\mu}(x,t)} w(\phi(x,t))\bigg\}
\label{O3prop3b} \\
& =  \int_0^t e^{-\alpha s-\Lambda^{\hat{\mu}(x)}(x,s)}\biggl[ -\rho + f(\phi(x,s),\hat{\mu}(x,s))
 + \lambda Qh(\phi(x,s),\hat{\mu}(x,s))\biggr] ds \nonumber \\
& \phantom{=} + e^{-\alpha t-\Lambda^{\hat{\mu}(x)}(x,t)} w(\phi(x,t)),
\label{O4prop3b}
\end{align}
where $\hat{\Theta}(x) = (\hat{\mu}(x), \hat{\mu}_{\partial}(x))$.
\end{proposition}

\bigskip

The following theorem shows the existence of a feedback measurable selector for the
one-stage optimization problems associated to $\mathcal{T}_{\alpha}(\rho,h)$ and $\mathcal{R}_{\alpha}(\rho,h)$.
Its proof is presented in the appendix.
\begin{theorem}
\label{theo3main} Let $\alpha\geq 0$, $\rho \in \RR_{+}$ and $h\in\mathbb{M}(E)$ be bounded from below.
For $x\in E$ define
\begin{equation}
w(x)=\mathcal{R}_{\alpha}(\rho,h)(x), 
\label{O1prop3b}
\end{equation}
and suppose that $w(x)\in \RR$ for all $x\in E$.
Then $w\in\mathbb{M}^{ac}(E)$ and the feedback measurable selector $\widehat{u}_{\phi}(w,h) \in \mathcal{S}_{\mathbb{V}}$ (see item D2) of Definition \ref{defmes}) satisfies
the following one-stage optimization problems:
\begin{eqnarray}
\mathcal{R}_{\alpha}(\rho,h)(x) & = & \mathcal{T}_{\alpha}(\rho,h)(x)\nonumber \\
& = & -\rho \mathcal{L}_{\alpha}(x,\widehat{u}_{\phi}(w,h)(x)) +
L_{\alpha}f(x,\widehat{u}_{\phi}(w,h)(x))+H_{\alpha}r(x,\widehat{u}_{\phi}(w,h)(x))\nonumber \\
& & +G_{\alpha}h(x,\widehat{u}_{\phi}(w,h)(x)).
 \label{minpro}
\end{eqnarray}
\end{theorem}

\section{Optimality equation for the discounted case}
\label{aux2}
In this section we consider the discounted optimal control problem (\ref{discostvf}) and, under the assumptions made in the previous
sections, we derive an optimality equation for this problem.
As usual in this kind of problem we characterize first the optimality equation for the truncated on the jump times $T_m$
problems (\ref{discostm}) and then take the limit as $m\rightarrow \infty$.

\bigskip
\noindent Throughout this section we consider $\alpha>0$ fixed. For any $g\in \mathbb{M}(E)^{+}$, we set $\mathcal{W}g$ as the
function on $E$ defined as
\begin{equation}
\mathcal{W}g(x) = \mathcal{R}_\alpha(0,g)(x), \label{OpT}
\end{equation}
for $x\in E$. The following proposition is an immediate consequence of the results derived in the previous section.

\begin{proposition}
\label{Tg}
For $g\in \mathbb{M}(E)^{+}$ consider $w=\mathcal{W}g$ and suppose that for all $x\in E$, $w(x)\in\RR$. Then $w\in \mathbb{M}(E)^{+}$ and
$\widehat{u}_{\phi}(w,g)\in \mathcal{S}_{\mathbb{V}}$ (see item D2) of Definition \ref{defmes}) satisfies
\begin{align}
w(x)  & =  L_{\alpha}f(x,\widehat{u}_{\phi}(w,g)(x))+H_{\alpha}r(x,\widehat{u}_{\phi}(w,g)(x))
+G_{\alpha}g(x,\widehat{u}_{\phi}(w,g)(x)).
\label{Tg1}
\end{align}
\end{proposition}
\textbf{Proof:} From Theorem \ref{theo3main}, we obtain the first equality and that $\widehat{u}_{\phi}(w,g)$ satisfies equation (\ref{Tg1}). Now
applying Proposition \ref{mes3} and by using the fact that $f\in \mathbb{M}(\widebar{E}\times\mathbb{U})^{+}$ and $r\in
\mathbb{M}(\partial E\times\mathbb{U})^{+}$, we obtain that $\mathcal{W}g\in \mathbb{M}(E)^{+}$. \hfill $\Box$

\bigskip

\noindent Define the sequence of functions $(v_m)_{m\in \NN}$ as
\begin{equation}
v_{m+1} = \mathcal{W}v_m,\,\,\,\, v_0=0, \label{v_m}
\end{equation}
We have the following proposition.

\begin{proposition}
\label{propvm} For all $x\in E$ and $m\in \NN$ we have that $\ds v_m(x) = \inf_{U\in\mathcal{U}}\mathcal{D}^{\alpha}_{m}(U,x)$.
\end{proposition}
\noindent {\textbf{Proof}:} It follows from the same lines as the proof of the Propositions \ref{prop2a} and \ref{prop2b}. \hfill
$\Box$

\bigskip

\noindent Since $\ds v_m(x) = \inf_{U\in\mathcal{U}}\mathcal{D}^{\alpha}_{m}(U,x)\leq \mathcal{J}_{\mathcal{D}}^{\alpha}(x)$,
the functions $v_m\in \mathbb{M}(E)^{+}$ and are non-decreasing.
Consequently, there exists $v\in \mathbb{M}(E)^{+}$ such that $v_m \uparrow v$, and it follows that $v\leq \mathcal{J}_{\mathcal{D}}^{\alpha}$.
We need the following propositions:

\begin{proposition}
\label{propge} If $h\in \mathbb{M}(E)^{+}$ is such that $h(x) \geq
\mathcal{W}h(x)$ then $h(x)\geq
\mathcal{J}_{\mathcal{D}}^{\alpha}(x)$.
\end{proposition}
\noindent {\textbf{Proof}:} By using Theorem \ref{theo3main} with $\rho=0$, we obtain that there exists $\widehat{u}_{\phi}\in
\mathcal{S}_{\mathbb{V}}$ such that
\begin{align*}
h(x) & \geq \mathcal{T}_\alpha(0,g)(x) =
L_{\alpha}f(x,\widehat{u}_{\phi}(x))+H_{\alpha}r(x,\widehat{u}_{\phi}(x))+G_{\alpha}h(x,\widehat{u}_{\phi}(x)).
\end{align*}
Define $w(x)=
L_{\alpha}f(x,\widehat{u}_{\phi}(x))+H_{\alpha}r(x,\widehat{u}_{\phi}(x))+G_{\alpha}h(x,\widehat{u}_{\phi}(x))$.
Clearly $w(x)\geq 0$. Moreover, the hypotheses of Proposition \ref{prop2b} are satisfied with $\rho=0$. Consequently, it follows
that there exists $\widehat{U} \in \mathcal{U}$ such that for all $m\in \NN$, $(t,x,k)\in \RR_{+}\times E\times \NN$,
\begin{align*}
E^{\widehat{U}}_{(x,k)} \biggl[  \int_{0}^{t\wedge T_{m}}e^{-\alpha s} & \Bigl[ f\bigl(X(s),\hat{u}(N(s),Z(s),\tau(s)) \bigr) \Bigr] ds \nonumber \\
& + \int_{0}^{t\wedge T_{m}} e^{-\alpha s} r\bigl(X(s-),\hat{u}_{\partial}(N(s),X(s-)) \bigr) dp^{*}(s) \biggr] \nonumber \\
& \leq h(x).
\end{align*}
From assumption \ref{Hypjump} (which implies that $T_{m}\rightarrow \infty$ $P^{\widehat{U}}$ a.s.), we have that
\begin{align*}
E^{\widehat{U}}_{(x,k)} \biggl[ & \int_{0}^{t} e^{-\alpha s} \Bigl[ f\bigl(X(s),\hat{u}(N(s),Z(s),\tau(s)) \bigr) \Bigr] ds \nonumber \\
& + \int_{0}^{t}e^{-\alpha s} r\bigl(X(s-),\hat{u}_{\partial}(N(s),X(s-)) \bigr) dp^{*}(s) \biggr] \nonumber \\
& \leq  h(x),
\end{align*}
and taking the limit as $t\rightarrow \infty$ we obtain that $h(x)\geq \mathcal{J}_{\mathcal{D}}^{\alpha}(x)$. \hfill $\Box$

\begin{proposition}
\label{propig} We have that $v(x) = \mathcal{W}v(x)$.
\end{proposition}
\noindent {\textbf{Proof}:} Let us show first that $v(x) \leq \mathcal{W}v(x)$. By using the definition of $\mathcal{W}$ we have
for any $\Upsilon\in \mathbb{V}^{r}(x)$ that \begin{equation*}
v_{m+1}(x) \leq
L_{\alpha}f(x,\Upsilon)+H_{\alpha}r(x,\Upsilon)+G_{\alpha}v_{m}(x,\Upsilon).
\end{equation*}
Taking the limit as $m\uparrow\infty$ and from the monotone convergence theorem we get that
\begin{align*}
v(x) & = \lim_{m\rightarrow \infty} v_{m+1}(x) \leq L_{\alpha}f(x,\Upsilon)+H_{\alpha}r(x,\Upsilon)+ \lim_{m\rightarrow \infty} G_{\alpha} v_{m}(x,\Upsilon)\\
&=  L_{\alpha}f(x,\Upsilon)+H_{\alpha}r(x,\Upsilon)+ G_{\alpha}
v(x,\Upsilon)
\end{align*}
showing that $v(x) \leq \mathcal{W}v(x)$. From Proposition \ref{Tg}, there exists for any $m\in \NN$, $u_{\phi}^m \in
\mathcal{S}_{\mathbb{V}}$ such that
\begin{equation}
\mathcal{W}v_m(x) =
L_{\alpha}f(x,u_{\phi}^{m}(x))+H_{\alpha}r(x,u_{\phi}^{m}(x))+G_{\alpha}v_m(x,u_{\phi}^{m}(x)).
\label{eqconv}
\end{equation}
Fix $x\in E$. 
Since $u_{\phi}^{m}(x)\in \mathbb{V}(x)\subset \mathbb{V}^{r}(x)$
and $\mathbb{V}^{r}(x)$ is compact we can find a further
subsequence, still written as $u_{\phi}^{m}(x)$ for notational
simplicity, such that $u_{\phi}^{m}(x) \rightarrow \hat{\Theta}
\in \mathbb{V}^{r}(x)$. From Proposition \ref{lemlsc1},
\begin{align}\label{selehadis}
v(x) &= \lim_{m\rightarrow \infty} v_{m+1}(x)\nonumber\\& =
\liminf_{m\rightarrow
\infty}\Bigl\{L_{\alpha}f(x,u_{\phi}^{m}(x))+H_{\alpha}r(x,u_{\phi}^{m}(x))+G_{\alpha}v_{m}(x,u_{\phi}^{m}(x))\Bigr\}
\nonumber\\& \geq
L_{\alpha}f(x,\hat{\Theta})+H_{\alpha}r(x,\hat{\Theta})+G_{\alpha}v(x,\hat{\Theta})
\geq \mathcal{R}_\alpha(0,g)(x) = \mathcal{W}v(x),
\end{align}
giving the result. \hfill $\Box$

\bigskip

\noindent Finally we have the following theorem characterizing the optimality equation for the discounted optimal control problem (\ref{discostvf}) and showing the convergence of the truncated problems.
\begin{theorem}
\label{theodis} We have that $v_n\uparrow
\mathcal{J}_{\mathcal{D}}^{\alpha}$ and
$\mathcal{J}_{\mathcal{D}}^{\alpha}(x) =
\mathcal{W}\mathcal{J}_{\mathcal{D}}^{\alpha}(x)$.
\end{theorem}
\noindent {\textbf{Proof}:} All we need to show is that
$\mathcal{J}_{\mathcal{D}}^{\alpha}(x) \leq v(x)$. But this is
immediate from Propositions \ref{propig} and \ref{propge}. \hfill
$\Box$

\section{The vanishing approach}\label{main3}
\noindent
In general it is hard to obtain a solution for the discrete-time optimality equation (see equation (\ref{OptEq0})). A common approach is to deal with an optimality inequality of the kind $h\geq \mathcal{T}(\rho,h)$.
We present sufficient conditions for the existence of a solution for this inequality, using the so-called vanishing discount approach (see Theorem \ref{main}).
Combining this result with the connection between the integro-differential equation and the discrete-time equation
we obtain our final main result that shows the existence of an ordinary optimal feedback control for the long run average cost (see Theorem \ref{main}).
First we have the following result, which traces a parallel with the Abelian Theorem (see \cite{hernandez96}).

\begin{proposition}
\label{abelian}
We have that
$\limsup_{\alpha \downarrow 0} \alpha \mathcal{J}_{\mathcal{D}}^{\alpha}(x) \leq \mathcal{J}_{\mathcal{A}}(x)$.
\end{proposition}
\noindent {\textbf{Proof}:}
See Theorem 1, chapter 5 in \cite{widder41}.
\hfill $\Box$

\bigskip

\noindent
We shall add the following assumptions for the discounted problems:

\begin{assumption}
\label{bounded}
There exists a state $x_0\in E$, numbers $\beta>0$, $C \geq 0$, $K_h \geq 0$, and a nonnegative function $b(.)$  such that for all $x\in E$ and $\alpha \in (0,\beta]$,
$\rho_{\alpha} \leq C,$
where $\rho_{\alpha} = \alpha \mathcal{J}_{\mathcal{D}}^{\alpha}(x_0)$ and
$-K_h \leq h_{\alpha}(x) \leq b(x)$
where $h_\alpha(x) = \mathcal{J}_{\mathcal{D}}^{\alpha}(x)-\mathcal{J}_{\mathcal{D}}^{\alpha}(x_0)$.
\end{assumption}

\bigskip

\noindent
We have the following propositions:
\begin{proposition}
\label{lemLas}
There exists a decreasing sequence of positive numbers $\alpha_k\downarrow 0$ such that $\rho_{\alpha_k}\rightarrow \rho$ and for all $x\in E$,
$\lim_{k\rightarrow \infty} \alpha_k  \mathcal{J}^{\alpha_k}(x) = \rho.$
\end{proposition}
\textbf{Proof:} See Lemma in \cite{hernandez96}, page 88.
\hfill $\Box$

\begin{proposition}
\label{lemlsc2}
Set $\ds h = \liminf_{k\rightarrow \infty} h_{\alpha_k}$. Then for all $x\in E$, $h(x) \geq -K_h$ and $h(x)\geq\mathcal{T}(\rho,h)(x)$.
\end{proposition}
\textbf{Proof:}
From Proposition \ref{Tg} and Theorem \ref{theodis} we have that the following equation is satisfied for each $\alpha > 0$ and $x\in E$:
\begin{align}
h_\alpha(x) & =
\mathcal{T}_\alpha(\rho_\alpha,h_\alpha)(x)\nonumber\\&
= -\rho_{\alpha}\mathcal{L}_{\alpha}(x,u^{\alpha}_{\phi}(x))+ L_{\alpha}f(x,u^{\alpha}_{\phi}(x))+H_{\alpha}r(x,u^{\alpha}_{\phi}(x))+G_{\alpha}h_{\alpha}(x,u^{\alpha}_{\phi}(x)),
\label{eqhal}
\end{align}
for $u^{\alpha}_{\phi} \in \mathcal{S}_{\mathbb{V}}$.
\nl
For $x\in E$ fixed and for all $k\in \NN$, $u^{\alpha_{k}}_{\phi}(x)\in \mathbb{V}(x)\subset \mathbb{V}^{r}(x)$ and since $\mathbb{V}^{r}(x)$ is compact we can
find a further subsequence, still written as $u^{\alpha_{k}}_{\phi}(x)$ for notational simplicity, such that $u^{\alpha_{k}}_{\phi}(x) \rightarrow \hat{\Theta} \in \mathbb{V}^{r}(x)$.
Combining Proposition \ref{lemlsc1} and equation (\ref{eqhal}),
\begin{align}\label{seleha}
h(x) &= \liminf_{k\rightarrow \infty} h_{\alpha_k}(x)\nonumber\\&
= \liminf_{k\rightarrow \infty}\Bigl\{-\rho_{\alpha_k} \mathcal{L}_{\alpha_k}(x,u^{\alpha_{k}}_{\phi}(x)) + L_{\alpha_k}f(x,u^{\alpha_{k}}_{\phi}(x))
+H_{\alpha_k}r(x,u^{\alpha_{k}}_{\phi}(x))+G_{\alpha_k}h_{\alpha_k}(x,u^{\alpha_{k}}_{\phi}(x))\Bigr\} \nonumber\\
& \geq -\rho\mathcal{L}(x,\hat{\Theta}) + L f(x,\hat{\Theta})+Hr(x,\hat{\Theta})+Gh(x,\hat{\Theta}).
\end{align}
Therefore, from Theorem \ref{theo3main}, it follows that
\begin{align*}
h(x) & \geq \mathcal{R}(\rho,h)(x) = \mathcal{T}(\rho,h)(x)
\end{align*}
showing the result.
\hfill $\Box$

\bigskip

\noindent
Our final result establishes the existence of an optimal control strategy for the long run average cost problem. 

\bigskip
\noindent
Let $h$ and $\rho$ be as in Propositions \ref{lemLas} and \ref{lemlsc2}, and $w=\mathcal{T}(\rho,h)$.

\begin{theorem}
\label{main}
$\widehat{U}_\phi(w,h) \in\mathcal{U}$ as defined in D.3) of Definition \ref{defmes} is such that
\begin{eqnarray*}
\rho  = \mathcal{J}_{\mathcal{A}}(x) = \mathcal{A}(\widehat{U}_\phi(w,h),x).
\end{eqnarray*}
\end{theorem}
\textbf{Proof:}
Combining Theorem \ref{theo3main} and Proposition \ref{lemlsc2}, we obtain that the hypotheses of Proposition \ref{prop2b} are satisfied for $\alpha=0$.
Consequently, setting for simplicity $\widehat{U}=\widehat{U}_\phi(w,h)$, it follows that
\begin{align*}
E^{\widehat{U}}_{(x,k)} \biggl[ & \int_{0}^{t\wedge T_{m}} \Bigl[ f\bigl(X(s),\hat{u}(N(s),Z(s),\tau(s)) \bigr) \Bigr] ds
+ \int_{0}^{t\wedge T_{m}} r\bigl(X(s-),\hat{u}_{\partial}(N(s),X(s-)) \bigr) dp^{*}(s) \biggr] \nonumber \\
& + E^{\widehat{U}}_{(x,k)} \bigl[ w\bigl(X(t\wedge T_{m}) \bigr)\bigr] \leq  E^{\widehat{U}}_{(x,k)} \bigl[ \rho [t\wedge T_{m}]\bigr] +w(x).
\end{align*}
Combining Proposition \ref{lemlsc2} and assumption \ref{Hyp8a}, we obtain that $w(x)\geq -\rho K_{\xi} - K_h$.
Moreover, we have, from assumption \ref{Hypjump} that $T_{m} \rightarrow \infty$ $P^{\widehat{U}}$ a.s. .
Consequently,
\begin{align*}
E^{\widehat{U}}_{(x,k)} \biggl[ & \int_{0}^{t} \Bigl[ f\bigl(X(s),\hat{u}(N(s),Z(s),\tau(s)) \bigr) \Bigr] ds
+ \int_{0}^{t} r\bigl(X(s-),\hat{u}_{\partial}(N(s),X(s-)) \bigr) dp^{*}(s) \biggr] \nonumber \\
& \leq  \rho\, t + \rho K_\varrho + K_h + w(x),
\end{align*}
showing that $\rho \geq \mathcal{A}(\widehat{U},x)$. From Proposition \ref{abelian} and Proposition \ref{lemLas}, we have $\rho \leq \mathcal{J}_{\mathcal{A}}(x)$ completing the proof.
\hfill $\Box$


\section{Appendix}
In this appendix we present several technical results required
throughout the paper.

\subsection{Proofs of the results of section \ref{conrex}}

We start with the proof of Proposition \ref{mesura1}.

\noindent {\bf{Proof of Proposition \ref{mesura1}}:} Let $\{x_{n}\}_{n\in \NN}$ be a convergent sequence in $E$ with $\ds
\lim_{n\rightarrow \infty} x_{n} = x$, $\mu \in \mathcal{V}^{r}(x)$ and $\bigl(\widebar{\mu}_{n}\bigr)_{n\in \NN}$
be a sequence in $\mathcal{V}^{r}(x_{n})$. Define the sequence $\bigl(\mu_{n}\bigr)_{n\in \NN}$ in $\mathcal{V}^{r}$ by
\begin{eqnarray*}
\mu_{n}(t,F) & \doteq & \frac{1}{\alpha(n,t)}\Bigr[\frac{1}{n}\widebar{\mu}_{n}(t,F)+(1-\frac{1}{n})
\mu(t,F\inter \mathbb{U}(\phi(x_{n},t))) \Bigl],
\end{eqnarray*}
where $(t,F)\in \RR_{+}\times \mathcal{B}(\mathbb{U})$, and $\alpha(n,t)=\frac{1}{n}+(1-\frac{1}{n})
\mu(t,\mathbb{U}(\phi(x_{n},t)))$. \nl Clearly, we have that $\mu_{n} \in \mathcal{V}^{r}(x_{n})$. Moreover, by using the
hypothesis it is easy to check that
\begin{eqnarray*}
\lim_{n\rightarrow \infty} \int_{\RR_{+}} \int_{\mathbb{U}} g(t,u) \mu_{n}(t,du) dt = \int_{\RR_{+}} \int_{\mathbb{U}}  g(t,u)
\mu(t,du) dt,
\end{eqnarray*}
for all $g\in L^{1}(\RR_{+};\mathbb{C}(\mathbb{U}))$. Therefore, from Proposition D.2 in \cite{hernandez96} the multifunction $x\in E$
$\rightarrow$ $\mathcal{V}^{r}(x)\subset \mathcal{V}^{r}$ is lower semicontinuous. From Corollary III.3 in \cite{castaing77}, this
multifunction is measurable and so from Proposition D.4 in \cite{hernandez96} $\bigl\{ (x,\mu) : \mu \in \mathcal{V}^{r}(x),
x\in E \bigr\} \in \mathcal{B}(E)\times\mathcal{B}(\mathcal{V}^{r})$. Finally, by using assumption
\ref{Hyp2a} it can be shown easily that $\mathcal{K} \doteq \bigl\{ (x,\Theta) : \Theta \in \mathbb{V}^{r}(x), x\in E \bigr\}
\in \mathcal{B}(E)\times \mathcal{B}(\mathbb{V}^{r})$, showing the result. \hfill $\Box$

\bigskip
\noindent We present next the proof of Proposition \ref{mes3}. First we need the following lemma.

\begin{lemma}
\label{mes2} Let $h\in \mathbb{M}(\RR_{+}\times E\times \mathcal{V}^r \times \RR_{+} \times \mathbb{U})^{+}$. Then the
mapping
\begin{center}
$(t,x,\mu) \in \RR_{+}\times E\times \mathcal{V}^r$ $\rightarrow$ $\ds \int_{\RR_{+}} \int_{\mathbb{U}} h(t,x,\mu,s,u) \mu(s,du) ds
\in \RR_{+}$
\end{center}
is $\mathcal{B}(\RR_{+}\times E\times \mathcal{V}^r)$-measurable.
\end{lemma}
\noindent {\bf{Proof}:} Let $\mathcal{H}$ be the class of functions $h\in \mathbb{B}(\RR_{+}\times E\times \mathcal{V}^r
\times \RR_{+} \times \mathbb{U})$ such that the mapping
\begin{center}
$(t,x,\mu) \in \RR_{+}\times E\times \mathcal{V}^r$ $\rightarrow$ $\ds \int_{\RR_{+}} e^{-s} \int_{\mathbb{U}} h(t,x,\mu,s,u)
\mu(s,du) ds \in \RR$
\end{center}
is $\mathcal{B}(\RR_{+}\times E\times \mathcal{V}^r)$-measurable. This set is closed relative to addition, multiplication by
constants and bounded pointwise passage to the limit. Consider the class $\mathcal{C}$ of functions $h$ such that $h=h_{1} h_{2}$
with $h_{1}= \mathbb{B}(\RR_{+}\times E\times \mathcal{V}^r )$ and $h_{2}= \mathbb{B}(\RR_{+} \times \mathbb{U})$. Then from Lemma
A.3 in \cite{forwick04}, $\mathcal{C} \subset \mathcal{H}$.
Moreover, $\mathcal{C}$ is closed relative to multiplication. Applying Theorem T20 in \cite{meyer66}, it follows that
$\mathcal{H}$ contains $\mathbb{B}(\RR_{+}\times E\times \mathcal{V}^r \times \RR_{+} \times \mathbb{U})$. Consider $g\in
\mathbb{M}(\RR_{+}\times E\times \mathcal{V}^r \times \RR_{+}\times \mathbb{U})^{+}$ and define $g_{k}(t,x,\mu,s,u) =
\bigl(e^{s} g(t,x,\mu,s,u) \bigr) \wedge k$. Therefore the mapping
\begin{center}
$(t,x,\mu) \in \RR_{+}\times E\times \mathcal{V}^r$ $\rightarrow$ $\ds \int_{\RR_{+}}  e^{-s} \int_{\mathbb{U}}  g_{k}(t,x,\mu,s,u)
\mu(s,du) ds \in \RR_{+}$
\end{center}
is  $\mathcal{B}(\RR_{+}\times E\times \mathcal{V}^r)$-measurable. By using the monotone convergence theorem, it follows
\begin{eqnarray*}
\lim_{k\rightarrow \infty}\ds \int_{\RR_{+}}  e^{-s}\int_{\mathbb{U}}   g_{k}(t,x,\mu,s,u) \mu(s,du) dt =
\int_{\RR_{+}} \int_{\mathbb{U}} g(t,x,\mu,s,u) \mu(s,du) ds,
\end{eqnarray*}
showing the result. \hfill $\Box$

\bigskip

\noindent {\bf{Proof of Proposition \ref{mes3}}:} Applying Lemma \ref{mes2} to $h(t,x,\mu,s,u)=I_{\{s\leq t\}}
\lambda(\phi(x,s),u)$ implies that the mapping $\Lambda^{\mu}(x,t)$ defined on $\RR_{+}\times E\times
\mathcal{V}^r$ with value in $\RR$ is measurable with respect to $\mathcal{B}(\RR_{+}\times E\times \mathcal{V}^r)$. There is no
loss of generality in assuming that $g$ is positive. Clearly, $Qg(.,.)$ is measurable with respect to $\mathcal{B}( \widebar{E}
\times \mathbb{U})$. \nl Therefore, for the function $h(t,x,\mu,s,u)=I_{\{s\leq t_{*}(x)\}} e^{-\alpha s -
\Lambda^{\mu}(x,s)} \lambda(\phi(x,s),u) Qg(\phi(x,s),u)$, Lemma \ref{mes2} shows that the mapping
\begin{center}
$(x,\mu) \in E\times \mathcal{V}^r$ $\rightarrow$ $\ds
\int_0^{t_{*}(x)}e^{-\alpha s - \Lambda^{\mu}(x,s)}\lambda
Qg(\phi(x,s),\mu(s)) ds \in \RR_{+}$
\end{center}
is $\mathcal{B}(E\times \mathcal{V}^r)$-measurable. Moreover, the mapping
\begin{center}
$(x,\mu,u) \in E\times \mathcal{V}^r\times \mathbb{U}$
$\rightarrow$ $\ds e^{-\alpha t_{*}(x) -\Lambda^{\mu}(x,t_{*}(x))}
Qg(\phi(x,t_{*}(x)),u) \in \RR_{+}$
\end{center}
is clearly $\mathcal{B}(E\times \mathcal{V}^r\times \mathbb{U})$-measurable. Consequently, it follows that $G_{\alpha}g(x,\Theta)$ defined on $E\times \mathbb{V}^r$ with
value in $\RR$ is $\mathcal{B}(E\times \mathbb{V}^r)$-measurable. By using the same arguments, the same property can be shown for
the mappings $L_{\alpha}w_{1}(x,\Theta)$, and $H_{\alpha}w_{2}(x,\Theta)$. \hfill $\Box$

\subsection{Proofs of the results of section \ref{main1}}

The next lemma applies the semi-group property of the flow $\phi$ (see equation (\ref{pflow})) into the operators operators $G_{\alpha}$,
$L_{\alpha}$, and $H_{\alpha}$ (defined in equations (\ref{DefGr}), (\ref{DefLr}), and (\ref{DefHr})).
Recall also the definition of $[\Theta]_{t}$ in (\ref{DefUps-rt}).

\begin{lemma}
\label{lem2} For any $\alpha \geq 0$, $x\in E$, $t\in [0,t_{*}(x))$, $\Theta=\bigl( \mu, \mu_{\partial} \bigr)\in
\mathbb{V}^{r}(x)$, and $g\in \mathbb{M}(E)$ bounded from below,
we have that
\begin{eqnarray}
\mathcal{L}_\alpha(x,\Theta)  & = & \int_0^t e^{-\alpha s-\Lambda^{\mu}(x,s)} ds
+ e^{-\alpha t-\Lambda^{\mu}(x,t)} \mathcal{L}_\alpha(\phi(x,t),\bigl[\Theta\bigr]_{t}), \label{eqLcr} \\
L_\alpha f(x,\Theta) & = & \int_0^t e^{-\alpha s -\Lambda^{\mu}(x,s)} f(\phi(x,s),\mu(s))ds + e^{-\alpha t -\Lambda^{\mu}(x,t)} L_\alpha f(\phi(x,t),\bigl[\Theta\bigr]_{t}),\nonumber\\
H_\alpha r(x,\Theta) & = & e^{-\alpha t -\Lambda^{\mu}(x,t)} H_\alpha r(\phi(x,t),\bigl[\Theta\bigr]_{t}),\nonumber\\
G_{\alpha} g(x,\Theta) & = & \int_0^t e^{-\alpha s -\Lambda^{\mu}(x,s)} \lambda Qg(\phi(x,s),\mu(s)) ds 
+e^{-\alpha t -\Lambda^{\mu}(x,t)} G_{\alpha}g((\phi(x,t)),[\Theta]_{t}). \nonumber
\end{eqnarray}
\end{lemma}
\noindent {\bf{Proof}:} For any $x\in E$, and $\Theta=\bigl( \mu,
\mu_{\partial} \bigr)\in \mathbb{V}^{r}(x)$, by using the
semi-group property of $\phi$, we have for $t+s< t_{*}(x)$
\begin{eqnarray*}
\Lambda^{\mu}(x,t+s) & =  & \int_0^t \lambda(\phi(x,\theta),\mu(\ell))d\ell + \int_t^{t+s} \lambda(\phi(x,\theta),\mu(\ell))d\ell \nonumber \\
& = & \Lambda^{\mu}(x,t)   + \int_{0}^{s}
\lambda(\phi(\phi(x,t),\theta),\mu(\ell+t)) d\ell.
\end{eqnarray*}
Remark that $t_{*}(x)-t=t_{*}(\phi(x,t))$. Consequently, combining
the previous equation and Definition \ref{def1*}, we obtain for
$t\in [0,t_{*}(x))$
\begin{eqnarray*}
\mathcal{L}_\alpha (x,\Theta)  & = & \int_0^t e^{-\alpha s
-\Lambda^{\mu}(x,s)} ds
+ \int_{0}^{t_{*}(\phi(x,t))} e^{-\alpha (t+s) -\Lambda^{\mu}(x,t+s)}ds, \nonumber \\
& = & \int_0^t e^{-\alpha s -\Lambda^{\mu}(x,s)} ds
+ e^{-\alpha t-\Lambda^{\mu}(x,t)} \int_0^{t_{*}(\phi(x,t))} e^{-\alpha s-\int_{0}^{s} \lambda(\phi(\phi(x,t),\ell),\mu(\ell+t)) d\ell}ds, \nonumber \\
& = & \int_0^t e^{-\alpha s -\Lambda^{\mu}(x,s)} ds + e^{-\alpha t
-\Lambda^{\mu}(x,t)}
\mathcal{L}_\alpha(\phi(x,t),\bigl[\Theta\bigr]_{t}).
\end{eqnarray*}
showing equation (\ref{eqLcr}). The other equalities can be obtained by using similar arguments. \hfill $\Box$

\bigskip
\noindent We present next the proof of Proposition \ref{prop3a}.

\noindent {\bf{Proof of Proposition \ref{prop3a}}:}
From Lemma \ref{lem2}, it follows that for any $x\in E$, $t\in [0,t_{*}(x))$, and $\Upsilon=\bigl( \nu,
\nu_{\partial} \bigr)\in \mathbb{V}(x)$
\begin{align}
- & \rho \mathcal{L}_{\alpha}(x,\Upsilon)  +  L_{\alpha}f(x,\Upsilon) + H_{\alpha}r(x,\Upsilon) + G_{\alpha} h(x,\Upsilon) \nonumber \\
= & e^{-\alpha t-\Lambda^{\nu}(x,t)} \biggl[ -\rho \mathcal{L}_{\alpha}(\phi(x,t),\bigl[\Upsilon\bigr]_{t}) +
L_{\alpha}f(\phi(x,t),\bigl[\Upsilon\bigr]_{t}) + H_{\alpha}r(\phi(x,t),\bigl[\Upsilon\bigr]_{t})
+ G_{\alpha} h(\phi(x,t),[\Upsilon]_{t}) \biggr] \nonumber \\
& + \int_0^t e^{-\alpha s -\Lambda^{\nu}(x,s)}\biggl[ -\rho + f(\phi(x,s),\nu(s)) + \lambda(\phi(x,s),\nu(s))
Qh(\phi(x,s),\nu(s))\biggr] ds . \label{eq1prop3a}
\end{align}
Notice now that we must have
\begin{eqnarray}
w(\phi(x,t)) & = & -\rho \mathcal{L}_{\alpha}(\phi(x,t),\bigl[\hat{\Gamma}(x)\bigr]_{t})
+ L_{\alpha}f(\phi(x,t),\bigl[\hat{\Gamma}(x)\bigr]_{t}) + H_{\alpha}r(\phi(x,t),\bigl[\hat{\Gamma}(x)\bigr]_{t})  \nonumber \\
& &+ G_{\alpha} h(\phi(x,t),[\hat{\Gamma}(x)]_{t})
\label{eq2prop3a}
\end{eqnarray}
otherwise that would contradict the fact that the infimum is reached in equation (\ref{OpDr1a}) for $\hat{\Gamma}(x)$.
Consequently, by taking $\Upsilon = \hat{\Gamma}(x)$ in equation (\ref{eq1prop3a}) we obtain equation (\ref{eqOpDr1b}).
\nl
From assumption \ref{Hyp4a}, we have that for all $x\in E$ and $t\in [0,t_{*}(x))$, $e^{-\Lambda^{\hat{\gamma}(x)}(x,t)} >0$ and so
equation (\ref{eqOpDr1b}) implies that for all $x\in E$, $w(\phi(x,t))$ is absolutely continuous on $[0,t_{*}(x))$.
Since there exists a constant $K_{h}\in \RR_{+}$ such that $-K_{h}\leq h$, it is easy to obtain
$$\ds \int_0^{t_{*}(x)} e^{-\alpha s-\Lambda^{\hat{\gamma}(x)}(x,s)} [ -\rho + f(\phi(x,s),\hat{\gamma}(x,s)) + \lambda Qh(\phi(x,s),\hat{\gamma}(x,s)) ] ds \leq w(x)+K_{h}.$$
Consequently, if $t_{*}(x)<\infty$, by using assumption \ref{Hyp4a} the limit of $w(\phi(x,t))$ as $t\rightarrow t_{*}(x)$ exists in $\RR$, showing that
$w\in\mathbb{M}^{ac}(E)$.

\bigskip

\noindent Let $\nu\in \mathcal{V}(x)$. Define $\tilde{\nu}$ by $\tilde{\nu}(s)= I_{[0,t[}(s)\nu(s)+I_{[t,\infty[}(s)\hat{\gamma}(x,s)$.
Then $\tilde{\nu}\in \mathcal{V}(x)$ and $\Upsilon$ defined by $\bigl( \tilde{\nu}(x), \hat{\gamma}_{\partial}(x) \bigr)$ belongs to $\mathbb{V}(x)$ and satisfies
$\bigl[\Upsilon\bigr]_{t}=\bigl[\hat{\Gamma}\bigr]_{t}$.
Consequently, combining (\ref{OpDr1a}), (\ref{eq1prop3a}) and (\ref{eq2prop3a}), it follows that for all $\nu\in \mathcal{V}(x)$
\begin{align*}
w(x) & \leq  -\rho \mathcal{L}_{\alpha}(x,\Upsilon) +  L_{\alpha}f(x,\Upsilon) + H_{\alpha}r(x,\Upsilon) + G_{\alpha} h(x,\Upsilon) \nonumber\\
&  =  e^{-\alpha t-\Lambda^{\nu}(x,t)}  w(\phi(x,t)) + \int_0^t e^{-\alpha s-\Lambda^{\nu}(x,s)}\biggl[ -\rho + f(\phi(x,s),\nu(s)) \nonumber \\
& \phantom{=} + \lambda(\phi(x,s),\nu(s))
Qh(\phi(x,s),\nu(s))\biggr] ds .
\end{align*}
Now from the previous equation and  (\ref{eqOpDr1b}), we obtain equation (\ref{eqOpDr1a}). \hfill $\Box$

\bigskip
\noindent We present next the proof of Proposition \ref{prop2a}.

\noindent {\bf{Proof of Proposition \ref{prop2a}}:} Set $w=\mathcal{T}_{\alpha}(\rho,h)$. By hypothesis, $w$ is bounded from below, and so $J^{U}_{m}(t,x,k)$ is well defined. For $U=(u,u_{\partial})\in
\mathcal{U}$, defined for $\hat{y}=(x,z,s,n)\in \widehat{E}$, $\widehat{f}^{U}(\hat{y})=f(x,u(n,z,s))$,
$\widehat{r}^{U}(\hat{y})=r(x,u_{\partial}(n,z))$, $\widehat{h}(\hat{y})=h(x)$, $\widehat{w}(\hat{y})=w(x)$, and for
$t\in [0,t_{*}(x)]$ $\widehat{\Lambda}^{U}(y,t)=\Lambda^{U}(x,n,t)$. Clearly, we have that $J^{U}_{0}(t,x,k) =w(x) \geq h(x)$ for all $(t,x,k)\in
\RR_{+}\times E\times \NN$. Now assume that for $m\in \NN$, $J^{U}_{m}(t,x,k)  \geq h(x)$ for all $(t,x,k)\in \RR_{+}\times
E\times \NN$. Define $\hat{x}=(x,x,0,k)$, then
\begin{align*}
 J^{U}_{m+1}(t,x,k) &=  E^{U}_{(x,k)} \Biggl[  I_{\{t<T_{1}\}} \biggl( \int_{0}^{t}  e^{-\alpha s} \Bigl[ \widehat{f}^{U}(\widehat{\phi}(\hat{x},s))  -\rho \Bigr] ds + e^{-\alpha t} \widehat{w}(\widehat{\phi}(\hat{x},t)) \biggr) \nonumber \\
& + I_{\{t\geq T_{1}\}} \biggl( \int_{0}^{t\wedge T_{m+1}}
e^{-\alpha s} \Bigl[ \widehat{f}^{U}(\widehat{X}^{U}(s)) -\rho
\Bigr] ds
+ \int_{0}^{t\wedge T_{m+1}} e^{-\alpha s} \widehat{r}^{U}\bigl(\widehat{X}^{U}(s-) \bigr) dp^{*}(s) \nonumber \\
 &  + e^{-\alpha t\wedge T_{m+1}} \widehat{w}\bigl(\widehat{X}^{U}(t\wedge T_{m+1}) \bigr) \biggr) \Biggr].
\end{align*}
Therefore,
\begin{align}
 J^{U}_{m+1}(t,x,k) & =  E^{U}_{(x,k)} \Biggl[ I_{\{t<T_{1}\}} \biggl( \int_{0}^{t} e^{-\alpha s} \Bigl[ \widehat{f}^{U}(\widehat{\phi}(\hat{x},s))  -\rho \Bigr] ds + e^{-\alpha t}\widehat{w}(\widehat{\phi}(\hat{x},t)) \biggr) \nonumber \\
& + I_{\{t\geq T_{1}\}} \biggl( \int_{0}^{T_{1}} e^{-\alpha s}
\Bigl[ \widehat{f}^{U}(\widehat{\phi}(\hat{x},s)) -\rho \Bigr] ds
+ I_{\{T_{1}=t_{*}(x)\}} e^{-\alpha t_{*}(x) } \widehat{r}^{U}\bigl(\widehat{\phi}(\hat{x},t_{*}(x))\bigr)  \biggr) \nonumber \\
& + I_{\{t\geq T_{1}\}} \biggl( \int_{T_{1}}^{t\wedge T_{m+1}}
e^{-\alpha s}\Bigl[ \widehat{f}^{U}(\widehat{X}^{U}(s)) -\rho
\Bigr] ds
+ \int_{T_{1}}^{t\wedge T_{m+1}} e^{-\alpha s} \widehat{r}^{U}\bigl(\widehat{X}^{U}(s-) \bigr) dp^{*}(s) \nonumber \\
 &  +  e^{-\alpha t\wedge T_{m+1}} \widehat{w}\bigl(\widehat{X}^{U}(t\wedge T_{m+1}) \bigr) \biggr) \Biggr].
 \label{eq1prop2a}
\end{align}
However, by using the strong Markov property of the process
$\{\widehat{X}^{U}(t)\}$, it follows that
\begin{align}
I_{\{t\geq T_{1}\}} & e^{-\alpha T_1} J^{U}_{m}(t-T_{1},\widehat{X}^{U}_{1},k+1) = E^{U}_{(x,k)}\Biggl[   I_{\{t\geq
T_{1}\}} \biggl(  \int_{T_{1}}^{t\wedge T_{m+1}}e^{-\alpha s}\Bigl[ \widehat{f}^{U}(\widehat{X}^{U}(s)) -\rho \Bigr] ds \nonumber \\
& + \int_{T_{1}}^{t\wedge T_{m+1}} e^{-\alpha s}\widehat{r}^{U}\bigl(\widehat{X}^{U}(s-) \bigr) dp^{*}(s) +
e^{-\alpha t\wedge T_{m+1}}\widehat{w}\bigl(\widehat{X}^{U}(t\wedge T_{m+1}) \bigr) \biggr) |
\mathcal{F}^{\widehat{X}^{U}}_{T_{1}}\Biggr]. \label{eq2prop2a}
\end{align}
Combining equations (\ref{eq1prop2a}) and (\ref{eq2prop2a}), and the fact that $I_{\{t\geq T_{1}\}} J^{U}_{m}(t-T_{1},\widehat{X}^{U}_{1},k+1)\geq
I_{\{t\geq T_{1}\}} \widehat{h}(\widehat{X}^{U}_{1})$ we obtain
\begin{align}
 J^{U}_{m+1}(t,x,k) \geq  E^{U}_{(x,k)} \Biggl[ & \int_{0}^{t\wedge T_{1}} e^{-\alpha s} \Bigl[ \widehat{f}^{U}(\widehat{\phi}(\hat{x},s))  -\rho \Bigr] ds + I_{\{t<T_{1}\}} e^{-\alpha t} \widehat{w}(\widehat{\phi}(\hat{x},t))  \nonumber \\
& + I_{\{t\geq T_{1} =t_{*}(x) \}} e^{-\alpha t_{*}(x)}\widehat{r}^{U}\bigl(\widehat{\phi}(\hat{x},t_{*}(x))\bigr) +
I_{\{t\geq T_{1}\}}  e^{-\alpha T_1}\widehat{h}(\widehat{X}^{U}_{1}) \Biggr]. \label{eq3prop2a}
\end{align}
However,
\begin{align}
E^{U}_{(x,k)} \Biggl[  &\int_{0}^{t\wedge T_{1}}  e^{-\alpha s}\Bigl[ \widehat{f}^{U}(\widehat{\phi}(\hat{x},s))  -\rho \Bigr] ds
+ I_{\{t<T_{1}\}} e^{-\alpha t}\widehat{w}(\widehat{\phi}(\hat{x},t))  \Biggr] \nonumber \\ & =
\int_{0}^{t\wedge t_{*}(x)} \Bigl[\widehat{f}^{U}(\widehat{\phi}(\hat{x},s))  -\rho \Bigr]
e^{-\alpha s-\widehat{\Lambda}^{U}(\hat{x},s)} ds + I_{\{t<t_{*}(x)\}}e^{-\alpha t -\widehat{\Lambda}^{U}(\hat{x},t)}
\widehat{w}(\widehat{\phi}(\hat{x},t)),
\end{align}
and
\begin{align}
& E^{U}_{(x,k)} \Biggl[ I_{\{t\geq T_{1}\}}  e^{-\alpha T_1}\widehat{h}(\widehat{X}^{U}_{1})
+ I_{\{t\geq T_{1} =t_{*}(x) \}} e^{-\alpha t_{*}(x)} \widehat{r}^{U}\bigl(\widehat{\phi}(\hat{x},t_{*}(x))\bigr) \Biggr] \nonumber \\
& = e^{-\alpha t_{*}(x) - \widehat{\Lambda}^{U}(\hat{x},t_{*}(x))}\widehat{r}(\widehat{\phi}(\hat{x},t_{*}(x))) I_{\{t\geq t_{*}(x)
\}} +  \int_{0}^{t\wedge t_{*}(x)} \widehat{Q}^{U}\widehat{h}(\widehat{\phi}(\hat{x},s))
\widehat{\lambda}^{U}(\widehat{\phi}(\hat{x},s)) e^{-\alpha s-\widehat{\Lambda}^{U}(\hat{x},s)} ds
\nonumber \\
&\phantom{=} + e^{-\alpha t_{*}(x)-\widehat{\Lambda}^{U}(\hat{x},t_{*}(x))} \widehat{Q}^{U}
\widehat{h}(\widehat{\phi}(\hat{x},t_{*}(x))) I_{\{t\geq t_{*}(x)\}}. \label{eq4prop2a}
\end{align}
Combining equations (\ref{eq3prop2a})-(\ref{eq4prop2a}), it follows that for $t\in \RR_{+}$
\begin{align}
 J^{U}_{m+1}&(t,x,k) \geq  \int_{0}^{t\wedge t_{*}(x)} \Bigl[ \widehat{f}^{U}(\widehat{\phi}(\hat{x},s))  -\rho
 + \widehat{Q}^{U} \widehat{h}(\widehat{\phi}(\hat{x},s)) \widehat{\lambda}^{U}(\widehat{\phi}(\hat{x},s))\Bigr] e^{-\alpha s -\widehat{\Lambda}^{U}(\hat{x},s)} ds \nonumber \\
& \phantom{\leq}  + I_{\{t\geq t_{*(x)}\}} e^{-\alpha t_{*}(x)-\widehat{\Lambda}^{U}(\hat{x},t_{*}(x))} \Bigl[
\widehat{Q}^{U}\widehat{h}(\widehat{\phi}(\hat{x},t_{*}(x)))+ \widehat{r}(\widehat{\phi}(\hat{x},t_{*}(x))) \Bigr] \nonumber \\
& \phantom{\leq} + I_{\{t<t_{*(x)}\}} e^{-\alpha t -\widehat{\Lambda}^{U}(\hat{x},t)} \widehat{w}(\widehat{\phi}(\hat{x},t)) \nonumber \\
& =\int_{0}^{t\wedge t_{*}(x)} e^{-\alpha s -\Lambda^{\nu_{k}}(x,s)}\biggl[ -\rho + f(\phi(x,s),\nu_{k}(s)) + \lambda(\phi(x,s),\nu_{k}(s)) Qh(\phi(x,s),\nu_{k}(s))\biggr] ds \nonumber \\
& \phantom{=} + I_{\{t\geq t_{*(x)}\}} e^{-\alpha t_{*}(x) -\Lambda^{\nu_{k}}(x,t_{*}(x))} \Bigl[ Qh(\phi(x,t_{*}(x)),u_{\partial}(k,x)) +r(\phi(x,t_{*}(x)),u_{\partial}(k,x)) \Bigr] \nonumber \\
& \phantom{=}  +  I_{\{t<t_{*(x)}\}} e^{-\alpha t-\Lambda^{\nu_{k}}(x,t)}w(\phi(x,t))  , \label{eq5prop2a}
\end{align}
with $\nu_{k}(.)=u(k,x,.)$. Clearly, $\nu_{k}(.)\in \mathcal{V}(x)$. Now if $t<t_{*}(x)$, then by applying Proposition
\ref{prop3a}, it follows that $J^{U}_{m+1}(t,x,k)\geq w(x)\geq h(x)$. If $t\geq t_{*}(x)$, then by defining $\Upsilon =
(\nu_{k},u_{\partial}(k,x))\in \mathbb{V}(x)$ and by using equation (\ref{OpD3a}) we have
\begin{align*}
 J^{U}_{m+1}(t,x,k) & \geq \int_{0}^{t_{*}(x)} e^{-\alpha s -\Lambda^{\nu_{k}}(x,s)}\biggl[ -\rho + f(\phi(x,s),\nu_{k}(s)) + \lambda(\phi(x,s),\nu_{k}(s)) Qh(\phi(x,s),\nu_{k}(s))\biggr] ds \nonumber \\
& \phantom{=} + e^{-\alpha t_{*}(x) -\Lambda^{\nu_{k}}(x,t_{*}(x))} \Bigl[ Qh(\phi(x,t_{*}(x)),u_{\partial}(k,x)) +r(\phi(x,t_{*}(x)),u_{\partial}(k,x)) \Bigr] \nonumber \\
& = -\rho \mathcal{L}_{\alpha}(x,\Upsilon) + L_{\alpha}f(x,\Upsilon)+H_{\alpha}r(x,\Upsilon)+G_{\alpha}
h(x,\Upsilon) \geq w(x)\geq h(x),
\end{align*}
showing the result. \hfill $\Box$

\bigskip
\noindent We present next the proof of Proposition \ref{prop2b}.

\noindent {\bf{Proof of Proposition \ref{prop2b}}:} From Lemma A.3 in \cite{forwick04} and item $(ii)$ of Lemma 3 in
\cite{yushkevich80} it follows that the mapping $\hat{u}$ defined by $\hat{u}(n,x,t)=\hat{\gamma}(x,t)$ belongs $\mathbb{M}(\NN
\times E\times \RR_{+};\mathbb{U})$ since $\hat{\Gamma}=(\hat{\gamma}, \hat{\gamma}_{\partial}) \in
\mathcal{S}_{\mathbb{V}}$. Clearly, $u(n,x,t)\in \mathbb{U}(\phi(x,t))$. Moreover, $\hat{u}_{\partial}$ defined by
$\hat{u}_{\partial}(n,x)=\hat{\gamma}_{\partial}(x)$ belongs $\mathbb{M}(\NN \times E;\mathbb{U})$ and satisfies
$u_{\partial}(n,x)\in \mathbb{U}(\phi(x,t_{*}(x)))$. Therefore, $\widehat{U}=(\hat{u},\hat{u}_{\partial})$ belongs to
$\mathcal{U}$.

\bigskip

\noindent Set $w=\mathcal{T}_{\alpha}(\rho,h)$.
It is easy to check that $w$ is bounded from below, and so $J^{\widehat{U}}_{m}(t,x,k)$ is well defined.
Remark that $J^{\widehat{U}}_{0}(t,x,k)=w(x)\leq h(x)$ for all $(t,x,k)\in \RR_{+}\times E\times \NN$. Now assume that
for $m\in \NN$, $J^{\widehat{U}}_{m}(t,x,k) \leq h(x)$ for all $(t,x,k)\in \RR_{+}\times E\times \NN$, then by using this
hypothesis, it follows that the inequalities in equations (\ref{eq3prop2a}), (\ref{eq5prop2a}) can be inverted for the
control process given by $\widehat{U}$. Consequently, if $t<t_{*}(x)$, the last statement of Proposition \ref{prop3a}
implies $J^{\widehat{U}}_{m+1}(t,x,k)\leq w(x) \leq  h(x)$. If $t\geq t_{*}(x)$, then
\begin{eqnarray*}
 J^{\widehat{U}}_{m+1}(t,x,k) \leq -\rho \mathcal{L}_\alpha(x,\hat{\Gamma}(x)) + L_\alpha f(x,\hat{\Gamma}(x))+H_\alpha r(x,\hat{\Upsilon}(x))+G_{\alpha}h(x,\hat{\Gamma}(x))=w(x) \leq h(x),
\end{eqnarray*}
showing the desired result. \hfill $\Box$

\subsection{Proofs of the results of section \ref{aux1}}
\noindent We present next the proof of Proposition \ref{lemlsc1}.

\noindent \textbf{Proof of Proposition \ref{lemlsc1}:} \nl
\underline{Item a)} For all $0\leq t < t_{*}(x)$, if $t_{*}(x)= \infty$, and all $0\leq t \leq t_{*}(x)$, if $t_{*}(x)< \infty$,
we have from assumptions \ref{Hyp4a} and \ref{Hyp3a} that
\begin{eqnarray*}
\lim_{n\rightarrow \infty} \int_{0}^t \int_{\mathbb{U}(\phi(x,s))} \lambda(\phi(x,s),u) \mu_{n}(s,du) ds = \int_{0}^t
\int_{\mathbb{U}(\phi(x,s))}  \lambda(\phi(x,s),u) \mu(s,du) ds,
\end{eqnarray*}
or, in other words,
\begin{eqnarray*}
\lim_{n\rightarrow \infty} \Lambda^{\mu_n}(x,t) = \Lambda^{\mu}(x,t). 
\end{eqnarray*}
From items a) and b) of assumption \ref{Hyp8a}, we have $\ds e^{-\alpha_n t-\Lambda^{\mu_n}(x,t)} \leq
e^{-\int_0^{t}\xi(\phi(x,s))ds}$ and \nl $\ds \int_0^{t_{*}(x)}e^{-\int_0^t \xi(\phi(x,s))ds}dt < \infty$). Consequently, by
using the dominated convergence theorem we obtain
\begin{eqnarray*}
\lim_{n\rightarrow \infty} \mathcal{L}_{\alpha_n}(x,\Theta_n) = \int_0^{t_{*}(x)} \lim_{n\rightarrow \infty} e^{-\alpha_n
t-\Lambda^{\mu_n}(x,t)}dt = \int_0^{t_{*}(x)}e^{-\alpha t-\Lambda^{\mu}(x,t)}dt = \mathcal{L}_{\alpha}(x,\Theta),
\end{eqnarray*}
showing item $a)$.

\bigskip

\noindent \underline{Item b)} We have from assumption \ref{Hyp6bis} that there exists a sequence of measurable functions
$(f_j)_{j\in \NN}$ such that for all $y\in \widebar{E}$, $f_j(y,.)\uparrow f(y,.)$ and $f_j(y,.) \in
\mathbb{C}(\mathbb{U}(y))$. We have for $(n,j)\in \NN^{2}$, $x\in E$,
\begin{eqnarray*}
L_{\alpha_n} f_{j}(x,\Theta_n) & = & \int_0^{t_{*}(x)} \bigl[e^{-\alpha_{n} t -\Lambda^{\mu_{n}}(x,t)} - e^{-\alpha t -\Lambda^{\mu}(x,t)} \bigr] f_j(\phi(x,t),\mu_n(t))dt \\
&  & + \int_0^{t_{*}(x)} e^{-\alpha t -\Lambda^{\mu}(x,t)}f_j(\phi(x,t),\mu_n(s))dt.
\end{eqnarray*}
However, item $a)$ of assumption \ref{Hyp8a} gives
\begin{eqnarray*}
\bigl| e^{-\alpha_{n} t -\Lambda^{\mu_{n}}(x,t)} - e^{-\alpha t -\Lambda^{\mu}(x,t)} \bigr| f_j(\phi(x,t),\mu_n(t)) & \leq & 2
e^{-\int_0^{t}\xi(\phi(x,s))ds} \sup_{a\in \mathbb{U}(\phi(x,t))}f(\phi(x,t),a).
\end{eqnarray*}
By combining item $c)$ of assumption \ref{Hyp8a} and the dominated convergence theorem we obtain
 \begin{eqnarray*}
\lim_{n\rightarrow \infty}\int_0^{t_{*}(x)} \bigl[e^{-\alpha_{n} t -\Lambda^{\mu_{n}}(x,t)} - e^{-\alpha t -\Lambda^{\mu}(x,t)}
\bigr] f_j(\phi(x,t),\mu_n(t))dt = 0.
\end{eqnarray*}
Therefore,
\begin{eqnarray*}
 \lim_{n\rightarrow \infty} L_{\alpha_n} f_{j}(x,\Theta_n) & = & \int_0^{t_{*}(x)} e^{-\alpha t -\Lambda^{\mu}(x,t)} f_j(\phi(x,t),\mu(s))dt.
\end{eqnarray*}
However, remark that $\ds L_{\alpha_n} f(x,\Theta_n) \geq L_{\alpha_n} f_{j}(x,\Theta_n)$, and the result follows by using
the monotone convergence theorem.

\bigskip

\noindent \underline{Item c)} Let us consider first that $t_{*}(x)= \infty$. From item $b)$ of assumption\ref{Hyp8a} and remark \ref{vac}
\begin{eqnarray*}
e^{-\Lambda^{\mu_n}(x,t_{*}(x))} \leq e^{-\int_0^{t_{*}(x)}\xi(\phi(x,s))ds}=0, \mbox{ and }
e^{-\Lambda^{\mu}(x,t_{*}(x))} \leq e^{-\int_0^{t_{*}(x)}\xi(\phi(x,s))ds}=0.
\end{eqnarray*}
and the result follows immediately since $H_{\alpha_n}r(x,\Theta_n)= Hr(x,\Theta) = 0$. \nl Suppose now
that $t_{*}(x) < \infty$ and set $z=\phi(x,t_{*}(x))$. We have from assumption \ref{Hyp7bis} that there exists a sequence of
measurable functions $(r_j)_{j\in \NN}$ such that for all $y\in \partial E$, $r_j(y,.)\uparrow r(y,.)$ and $r_j(y,.) \in
\mathbb{C}(\mathbb{U}(y))$. Consequently, $r(z,\mu_{\partial,n}) \geq r_i(z,\mu_{\partial,n})$, and so $\ds \liminf_{n\rightarrow
\infty} r(z,\mu_{\partial,n}) \geq r_i(z,\mu_{\partial})$. From the monotone convergence theorem we obtain $c)$.

\bigskip

\noindent \underline{Item d)} First notice that
\begin{eqnarray}
\lim_{n\rightarrow \infty} L_{\alpha_n}\lambda(x,\Theta_n) = L\lambda(x,\Theta).
\label{eqcon2}
\end{eqnarray}
Indeed, notice that $L_{\alpha_n}(\lambda+\alpha_n)(x,\Theta_n) = 1 - e^{-\alpha_n t_{*}(x) -\Lambda^{\mu_n}(x,t_{*}(x))}$.
Considering first that $t_{*}(x)< \infty$, we have that
\begin{align*}
\lim_{n\rightarrow \infty} L_{\alpha_n} (\lambda)(x,\Theta_n) & =
\lim_{n\rightarrow \infty} L_{\alpha_n} (\lambda+\alpha_n)(x,\Theta_n) -\alpha \mathcal{L}_{\alpha}(x,\Theta)
\\& = 1 - e^{-\alpha t_{*}(x)-\Lambda^{\mu}(x,t_{*}(x))} -\alpha \mathcal{L}_{\alpha}(x,\Theta)
\\& = L (\lambda+\alpha) (x,\Theta)-\alpha \mathcal{L}_{\alpha}(x,\Theta) = L\lambda(x,\Theta).
\end{align*}
If $t_{*}(x)=\infty$ then  $1=L_{\alpha_n}\lambda(x,\Theta_n) = L \lambda(x,\Theta)$, showing (\ref{eqcon2}).
\nl
Set $\tilde{h}_{\alpha_k}= h_{\alpha_k}+K_h$, $\tilde{h}= h+K_h$ and $\tilde{g}_k = \inf_{j\geq k}\tilde{h}_{\alpha_j}$ (thus $\tilde{g}_k\uparrow \tilde{h}$ and $\tilde{g}_k\leq \tilde{h}_{\alpha_n}$ for $n\geq k$).
By hypothesis, $\tilde{g}_k(y) \geq 0$ for all $y\in E$.
We have that $\tilde{g}_k$ is the limit of a nondecreasing sequence of measurable bounded functions $\tilde{g}_{k,i} \in \mathbb{B}(E)$. Set $\lambda_m(y,a) = m\wedge\lambda(y,a)$.
From assumptions \ref{Hyp3a} and \ref{Hyp5a}, we have that for each $k$, $i$ and $m$ and $y\in E$, $\lambda_m Q\tilde{g}_{k,i}(y,.)$ is continuous on $\mathbb{U}(y)$.
Assumption \ref{Hyp8a} and the fact that for each $k,i$, $\tilde{g}_{k,i}$ is bounded above by, say $M_{k,i}$, yields
\begin{align*}
0 &\leq \int_0^{t_{*}(x)}  e^{-\int_0^t \xi(\phi(x,s))ds}\sup_{a\in \mathbb{U}(\phi(x,t))}(\lambda_m Q\tilde{g}_{k,i}(\phi(x,t),a))dt\\ &\leq
m \,M_{k,i} \int_0^{t_{*}(x)}  e^{-\int_0^t \xi(\phi(x,s))ds}dt < \infty.
\end{align*}
Since $(\lambda_m Q \tilde{g}_{k,i})(y,a)\geq 0$ and is continuous in $a$ we have from the proof of
b) (replacing $f$ by  $\lambda_m Q \tilde{g}_{k,i}$) that $\ds \liminf_{n\rightarrow \infty} L_{\alpha_n}(\lambda_m Q \tilde{g}_{k,i})(x,\Theta_n) \geq L(\lambda_m Q \tilde{g}_{k,i})(x,\Theta)$, and thus, recalling
that $\tilde{g}_{k,i}\leq \tilde{h}_{\alpha_n}$ for $n\geq k$ and $\lambda_m \leq \lambda$,
\begin{eqnarray*}
\liminf_{n\rightarrow \infty} L_{\alpha_n}(\lambda Q\tilde{h}_{\alpha_n})(x,\Theta_n) \geq L_{\alpha}(\lambda_m Q\tilde{g}_{k,i})(x,\Theta).
\end{eqnarray*}
From the monotone convergence theorem and taking the limit over $m,i,k$ we get that
\begin{eqnarray}
\liminf_{n\rightarrow \infty} L_{\alpha_n}(\lambda Q\tilde{h}_{\alpha_n})(x,\Theta_n) \geq L_{\alpha}(\lambda Q\tilde{h})(x,\Theta).
\label{eqcon3}
\end{eqnarray}
Notice now that
\begin{eqnarray*}
L_{\alpha_n}(\lambda Q\tilde{h}_{\alpha_n})(x,\Theta_n) =
L_{\alpha_n}(\lambda Q h_{\alpha_n})(x,\Theta_n) + K_h L_{\alpha_n}(\lambda)(x,\Theta_n)
\end{eqnarray*}
and similarly
\begin{eqnarray*}
L_{\alpha}(\lambda Q\tilde{h}_{\alpha})(x,\Theta) =
L_{\alpha} (\lambda Qh_{\alpha})(x,\Theta) + K_h L_{\alpha}(\lambda)(x,\Theta).
\end{eqnarray*}
By combining (\ref{eqcon2}) and (\ref{eqcon3}) we get that $\liminf_{n\rightarrow \infty} L_{\alpha_n}(\lambda Q h_{\alpha_n})(x,\Theta_n) \geq L_{\alpha}(\lambda Q h)(x,\Theta)$.
Using similar arguments as above and c) we can show that
\begin{eqnarray*}
\liminf_{n\rightarrow \infty} H_{\alpha_n}h_{\alpha_n}(x,\Theta_n) \geq  H h(x,\Theta).
\end{eqnarray*}
completing the proof of d).
\hfill $\Box$

\subsection{Proofs of the results of section \ref{main2}}

We present first the proof of Proposition \ref{prop3b}.
\bigskip

\noindent {\bf{Proof of Proposition \ref{prop3b}}:} From assumption \ref{Mesurability} and Proposition \ref{mes3}, it
follows that the mapping $V$ defined on $\mathcal{K}$ by
$$ V(x,\Theta)=-\rho \mathcal{L}_{\alpha}(x,\Theta) + L_{\alpha}f(x,\Theta)+H_{\alpha}r(x,\Theta)+G_{\alpha} h(x,\Theta) $$
is measurable. Moreover, by using Corollary \ref{corlsc1} it follows that for all $x\in E$, $V(x,.)$ is lower semicontinuous on
$\mathbb{V}^{r}(x)$. Recalling that $\mathbb{V}^{r}(x)$ is a compact subset of $\mathbb{V}^{r}$ and by using Proposition D.5 in
\cite{hernandez96}. we obtain that there exists $\hat{\Theta}\in \mathcal{S}_{\mathbb{V}^{r}}$ such equation (\ref{O2prop3b}) is
satisfied. \nl The rest of the proof is similar to the proof of Proposition \ref{prop3a} and it is, therefore, omitted. \hfill
$\Box$

\bigskip

\noindent Before presenting the proof of Theorem \ref{theo3main} we need the following auxiliary results.

\bigskip

\begin{lemma}
\label{Iexist} Assume that $w\in \mathbb{M}^{ac}(E)$. Then there exists a function $\mathcal{X}(w)$ in $\mathbb{M}(E)$ such that
for all $x\in E$, and $t\in [0,t_{*}(x))$
\begin{eqnarray}
w(\phi(x,t))-w(x) & = & \int_{0}^{t} \mathcal{X}(w)(\phi(x,s)) ds.
\label{derixt}
\end{eqnarray}
\end{lemma}
\noindent {\bf{Proof}:} Define
\begin{eqnarray*}
w^{+}(x) & \doteq & \limsup_{n\rightarrow +\infty} n\bigl[ w(\phi(x,t_{*}(x)\wedge \frac{1}{n+1})) -w(x)\bigr] \\
w^{-}(x) & \doteq & \liminf_{n\rightarrow +\infty} n\bigl[
w(\phi(x,t_{*}(x)\wedge \frac{1}{n+1})) -w(x)\bigr]
\end{eqnarray*}
Since $\bigl\{ w(\phi(x,t_{*}(x)\wedge \frac{1}{n+1}))
\bigr\}_{n\in\NN}$ is a sequence in $\mathbb{M}(E)$, then $w^{+}(x)$ and $w^{-}(x)$ are Borel measurable functions from $E$
into $\RR\cup\{\infty\}\cup\{-\infty\}$. Consequently, the set $\mathcal{D}_{w} \doteq \Bigl\{ x \in E : w^{+}(x)=w^{-}(x)
\Bigr\} \inter \Bigl\{ x \in E : w^{+}(x)\in \RR \Bigr\}$ belongs to  $\mathcal{B}(E)$.
 \nl
Define the function $\mathcal{X}(w)(x)$ by
\begin{eqnarray*}
\mathcal{X}(w)(x) = \left\{
\begin{array}{cr}
\frac{d w(\phi(x,t))}{dt}|_{t>0}, & \mbox{ if } x\in \mathcal{D}_{w}, \\
g(x), & \mbox{ otherwise, }
\end{array}
\right.
\end{eqnarray*}
where $g$ is any function in $\mathbb{M}(E)$. \nl Clearly $\mathcal{X}(w)$ belongs to $\mathbb{M}(E)$. Since $w\in
\mathbb{M}^{ac}(E)$, there exists a set $T^{w}_{x}\in \mathcal{B}([0,t_{*}(x)))$ such that $\mu_{leb}(
(T^{w}_{x})^{c}\cap[0,t_{*}(x)))=0$ and $w(\phi(.,x))$ admits derivatives in $T^{w}_{x}$. Consequently, for any $x\in E$, and
$t_{0}\in T^{w}_{x}$, we obtain that $\phi(x,t_{0})\in \mathcal{D}_{w}$, and
\begin{eqnarray*}
\lim_{\ds \mathop{\SS \epsilon \rightarrow 0}_{\SS \epsilon >0}} \frac{1}{\epsilon} \bigl[ w(\phi(x,t_{0}+\epsilon))
-w(\phi(x,t_{0}))\bigr] = \frac{dw(\phi(\phi(x,t_{0}),t))}{dt}|_{t>0} =
\mathcal{X}(w)(\phi(x,t_{0})).
\end{eqnarray*}
Therefore, $\mathcal{X}(w)$ satisfies (\ref{derixt})
showing the result. \hfill $\Box$

\bigskip

\begin{lemma}
\label{Aux1a} For any $\mu\in \mathcal{P}( \mathbb{U}(x))$ and $x\in \widebar{E}$, $\lambda(x,\mu)<\infty$.
\end{lemma}
\textbf{Proof:} From assumption \ref{Hyp1a}, $\mathbb{U}(x)$ is a compact subspace of $\mathbb{U}$, and from assumption
\ref{Hyp3a}, $\lambda(x,.):\mathbb{U}(x)\mapsto \RR_{+}$ is continuous. Therefore there exists $\hat{a}\in  \mathbb{U}(x)$ such that $\max_{a\in \mathbb{U}(x)}\lambda(x,a) = \lambda(x,\hat{a})$
and thus $0\leq \lambda(x,\mu) = \int_{\mathbb{U}(x)}\lambda(x,a) \mu(da) \leq  \lambda(x,\hat{a})$.
\hfill $\Box$

\bigskip

\begin{lemma}\label{Aux2a}
Suppose that $h\in \mathbb{M}(E)$ is bounded from below by $K_h$. Then
\begin{align}
\inf_{a\in \mathbb{U}(x)} \Bigl\{f(x,a) &+ \lambda(x,a)Q(h)(x,a) - \lambda(x,a)w(x) \Bigr\}\nonumber\\& = \inf_{\mu\in \mathcal{P}( \mathbb{U}(x)) }\Bigl\{f(x,\mu)-\lambda(x,\mu)w(x)+ \lambda Qh(x,\mu)\Bigr\}, \label{infmin2}\\
\inf_{a\in \mathbb{U}(\phi(x,t_{*}(x)))}\{ & r(\phi(x,t_{*}(x)),a)+Qh(\phi(x,t_{*}(x)),a)\} \nonumber\\& = \inf_{\mu\in
\mathcal{P}(\mathbb{U}(\phi(x,t_{*}(x))))}\{r(\phi(x,t_{*}(x)),\mu)+Qh(\phi(x,t_{*}(x)),\mu)\}.
\label{infmin1a}
\end{align}
\end{lemma}
\textbf{Proof:}
Set for simplicity, $\vartheta(x,a) = f(x,a) + \lambda(x,a)Q(h)(x,a) - \lambda(x,a)w(x)$.
Notice that from Lemma \ref{Aux1a}, for any $\mu \in \mathcal{P}( \mathbb{U}(x))$, $\lambda(x,\mu)< \infty$ and thus, recalling that $f$ and $h+K_h$ are positive,
\begin{align}
&f(x,\mu)-\lambda(x,\mu)w(x)+
\lambda Qh(x,\mu) = f(x,\mu)+ \lambda Q(h+K_h)(x,\mu) -\lambda(x,\mu)(w(x)+K_h)=\nonumber\\
& \int_{\mathbb{U}(x)}\left(f(x,a) + \lambda(x,a)Q(h+K_h)(x,a) - \lambda(x,a)(w(x)+K_h)\right)\mu(da) =
\int_{\mathbb{U}(x)}\vartheta(x,a) \mu(da). \label{eqaAux1a}
\end{align}
But as in Lemma 5.7 of \cite{forwick04}, we have that
\begin{align}
\inf_{a\in \mathbb{U}(x)} \vartheta(x,a) & = \inf_{\mu\in \mathcal{P}( \mathbb{U}(x))} \int_{\mathbb{U}(x)}\vartheta(x,a) \mu(da). \label{eqaAux1b}
\end{align}
Combining (\ref{eqaAux1a}) and (\ref{eqaAux1b}) we get (\ref{infmin2}). Similarly we have (\ref{infmin1a}).
\hfill $\Box$

\bigskip

\noindent We present next the proof of Theorem \ref{theo3main}.

\bigskip

\noindent \textbf{Proof of Theorem \ref{theo3main}:}
According to Proposition \ref{prop3b}, there exists $\hat{\Theta}\in \mathcal{S}_{\mathbb{V}^{r}}$ such that for all $x\in E$ and $t\in [0,t_{*}(x))$ we have
\begin{eqnarray}
e^{-\alpha t-\Lambda^{\hat{\mu}(x)}(x,t)}w(\phi(x,t))- w(x) & = & \int_0^t e^{-\alpha s -\Lambda^{\hat{\mu}(x)}(x,s)}\biggl[ \rho - f(\phi(x,s),\hat{\mu}(x,s)) \nonumber \\
& & - \lambda Qh(\phi(x,s),\hat{\mu}(x,s)) \biggr] ds,
\label{eq2theo3a}
\end{eqnarray}
where $\hat{\Theta}(x)=(\hat{\mu}(x), \hat{\mu}_{\partial}(x))$.
Since $w\in \mathbb{M}^{ac}(E)$, it follows from Lemma \ref{Iexist} that there exists a function $\mathcal{X}(w)$ in $\mathbb{M}(E)$
satisfying equation (\ref{derixt}).
Therefore,  we obtain from equation (\ref{eq2theo3a}) that
\begin{align*}
\mathcal{X}w(\phi(x,t))-[\alpha+\lambda(\phi(x,t),\hat{\mu}(x,t))] w(\phi(x,t))= -f(\phi(x,t),\hat{\mu}(x,t))
-\lambda Qh(\phi(x,t),\hat{\mu}(x,t))+\rho,
\end{align*}
$\eta-a.s.$ on $[0,t_{*}(x))$, implying that
\begin{align*}
-\mathcal{X} & w(\phi(x,t))+ \alpha w(\phi(x,t)) \nonumber \\
& \geq  \inf_{\mu \in \mathcal{P} \bigl(\mathbb{U}(\phi(x,t))\bigr)}  \Bigl\{f(\phi(x,t),\mu)-\lambda(\phi(x,t),\mu)w(\phi(x,t)) + \lambda Qh(\phi(x,t),\mu) \Bigr\} - \rho.
\end{align*}
However, remark that
\begin{align*}
\inf_{\mu \in \mathcal{P} \bigl( \mathbb{U}(\phi(x,t))\bigr)}  \Bigl\{f( & \phi(x,t),\mu)-\lambda(\phi(x,t),\mu)w(\phi(x,t))+ \lambda Qh(\phi(x,t),\mu) \Bigr\} - \rho \\
& =  \inf_{a\in \mathbb{U}(\phi(x,t))} \Bigl\{f(\phi(x,t),a)- \lambda(\phi(x,t),a) \bigl[ w(\phi(x,t)) - Qh(\phi(x,t),a) \bigr] \Bigr\} - \rho
\end{align*}
Consequently, by considering the measurable selector $\widebar{u}\in \mathcal{S}_{\mathbb{U}}$
given by $\widebar{u} = \widehat{u}(w,h)$ (see Definition \ref{defmes}, D1)), we have
\begin{align}
-\mathcal{X} & w(\phi(x,t))+ \alpha w(x) \nonumber \\
& \geq - \rho + f(\phi(x,t),\widebar{u}(\phi(x,t)))-\lambda(\phi(x,t),\widebar{u}(\phi(x,t))) \bigl[
w(\phi(x,t))-Qh(\phi(x,t),\widebar{u}(\phi(x,t))) \bigr], \label{eq3theo3a}
\end{align}
$\eta-a.s.$ on $[0,t_{*}(x))$.
Let $\Xi$ be the set in  $\mathcal{B}( [0,t_{*}(x)))$ such that the previous inequality is strict.
If $\eta(\Xi)>0$, then there would exist $t\in [0,t_{*}(x))$ such that
\begin{align*}
w(x) -e^{-(\alpha t+\overline{\Lambda}(x,t))} & w(\phi(x,t)) >
\int_0^t e^{-(\alpha s+\overline{\Lambda}(x,s))} \Bigl[ f(\phi(x,s),\widebar{u}(\phi(x,s))) \nonumber \\
& +\lambda(\phi(x,s),\widebar{u}(\phi(x,s))) Qh(\phi(x,s),\widebar{u}(\phi(x,s))) - \rho \Bigr] ds,
\end{align*}
where $\ds \overline{\Lambda}(x,t)$ denotes $\ds \int_{0}^{t}\lambda(\phi(x,s),\widebar{u}(\phi(x,s))) ds$.
However, this would lead to a contradiction with equation (\ref{O3prop3b}).
Thus we have
\begin{eqnarray*}
-\mathcal{X}w(\phi(x,t))+\alpha w(\phi(x,t)) & = & - \rho + f(\phi(x,t),\widebar{u}(\phi(x,t)))\nonumber \\
& & -\lambda(\phi(x,t),\widebar{u}(\phi(x,t))) \bigl[w(\phi(x,t))-Qh(\phi(x,t),\widebar{u}(\phi(x,t))) \bigr],
\end{eqnarray*}
$\eta-a.s.$ on $[0,t_{*}(x))$.
Consequently, for all $t\in [0,t_{*}(x))$ it follows that
\begin{align}
w(x) = e^{-(\alpha t+\overline{\Lambda}(x,t))} & w(\phi(x,t)) +
\int_0^t e^{-(\alpha s+\overline{\Lambda}(x,s))} \Bigl[ f(\phi(x,s),\widebar{u}(\phi(x,s))) \nonumber \\
& +\lambda(\phi(x,s),\widebar{u}(\phi(x,s))) Qh(\phi(x,s),\widebar{u}(\phi(x,s))) - \rho \Bigr] ds.
\label{espoir}
\end{align}

\noindent First consider the case in which $t_{*}(x)<\infty$. We obtain, by taking the limit as $t$ tends to $t_{*}(x)$ in the previous equation,
that the feedback measurable selector $\widehat{u}_{\phi}(w,h) \in \mathcal{S}_{\mathbb{V}}$ (see item D2) of Definition \ref{defmes}) satisfies:
\begin{align}
w(x) = & e^{-(\alpha t_{*}(x)+\overline{\Lambda}(x,t_{*}(x)))}  w(\phi(x,t_{*}(x))) 
-\rho \mathcal{L}_{\alpha}(x,\widehat{u}_{\phi}(w,h)(x)) + L_{\alpha}f(x,\widehat{u}_{\phi}(w,h)(x)) \nonumber \\
& +\int_0^{t_{*}(x)} e^{-(\alpha s+\overline{\Lambda}(x,s))} \lambda(\phi(x,s),\widebar{u}(\phi(x,s))) Qh(\phi(x,s),\widebar{u}(\phi(x,s)))ds.
\label{infmin1bis}
\end{align}
Define the control $\Theta(x)$ by $(\hat{\mu}(x), \mu)$ for $\mu \in
\mathcal{P}\bigl(\mathbb{U}(\phi(x,t_{*}(x)))\bigr)$. From equation (\ref{O1prop3b}), we obtain that
\begin{align*}
w(x) \leq &-\rho \mathcal{L}_{\alpha}(x,\Theta(x)) + L_{\alpha}f(x,\Theta(x))+H_{\alpha}r(x,\Theta(x))+G_{\alpha}
h(x,\Theta(x)),
\end{align*}
and by using the definition of $\Theta(x)$ and $\hat{\Theta}(x)$, it follows that
\begin{align}
w(x) \leq &-\rho \mathcal{L}_{\alpha}(x,\hat{\Theta}(x))+ L_{\alpha}f(x,\hat{\Theta}(x))
+\int_{0}^{t_{*}(x)} e^{-\alpha s -\Lambda^{\hat{\mu}(x)}(x,s)} \lambda Qh(\phi(x,s),\hat{\mu}(x,s)) ds \nonumber \\
& +e^{-\alpha t_{*}(x)-\Lambda^{\hat{\mu}(x)}(x,t_{*}(x))} \bigl[Qh(\phi(x,t_{*}(x)),\mu) + r(\phi(x,t_{*}(x)),\mu) \bigr].
\label{eq1coro2a}
\end{align}
From equation (\ref{O4prop3b}), we have that
\begin{align*}
w(x) = & \int_0^t e^{-\alpha s-\Lambda^{\hat{\mu}(x)}(x,s)}\biggl[ -\rho + f(\phi(x,s),\hat{\mu}(x,s))  + \lambda Qh(\phi(x,s),\hat{\mu}(x,s)) \biggr] ds \nonumber \\
& + e^{-\alpha t
-\Lambda^{\hat{\mu}(x)}(x,t)}w(\phi(x,t)).
\end{align*}
Since $w\in \mathbb{M}^{ac}(E)$, this yields that
\begin{align}
w(x) = & \lim_{t\rightarrow t_{*}(x)} \int_0^t e^{-\alpha s-\Lambda^{\hat{\mu}(x)}(x,s)}\biggl[ -\rho + f(\phi(x,s),\hat{\mu}(x,s))  + \lambda Qh(\phi(x,s),\hat{\mu}(x,s)) \biggr] ds \nonumber \\
& + \lim_{t\rightarrow t_{*}(x)} e^{-\alpha t -\Lambda^{\hat{\mu}(x)}(x,t)}w(\phi(x,t)) \nonumber \\
= &-\rho \mathcal{L}_{\alpha}(x,\hat{\Theta}(x)) +
L_{\alpha}f(x,\hat{\Theta}(x))
+\int_{0}^{t_{*}(x)} e^{-\alpha s -\Lambda^{\hat{\mu}(x)}(x,s)} \lambda Qh(\phi(x,s),\hat{\mu}(x,s)) ds \nonumber \\
& +e^{-\alpha t_{*}(x)-\Lambda^{\hat{\mu}(x)}(x,t_{*}(x))}w(\phi(x,t_{*}(x))).
\label{eq2coro2a}
\end{align}
From assumption \ref{Hyp4a}, we have that $e^{-\Lambda^{\hat{\mu}(x)}(x,t_{*}(x))} >0$. Therefore, combining
equations (\ref{eq1coro2a}), and (\ref{eq2coro2a}), it gives that for all $x\in E$, and $\mu \in
\mathcal{P}\bigl(\mathbb{U}(\phi(x,t_{*}(x)))\bigr)$
\begin{eqnarray*}
w(\phi(x,t_{*}(x))) & \leq  & Qh(\phi(x,t_{*}(x)),\mu)+r(\phi(x,t_{*}(x)),\mu).
\end{eqnarray*}
Clearly, by using equation (\ref{O2prop3b}), it can be claimed that the previous inequality becomes an equality for
$\mu=\hat{\mu}_{\partial}(x)$, implying that
\begin{eqnarray*}
w(\phi(x,t_{*}(x))) & = & \inf_{\mu\in \mathcal{P}(\mathbb{U}(\phi(x,t_{*}(x))))}\{r(\phi(x,t_{*}(x)),\mu)+Qh(\phi(x,t_{*}(x)),\mu)\} \\
& = & \inf_{a\in \mathbb{U}(\phi(x,t_{*}(x)))}\{  r(\phi(x,t_{*}(x)),a)+Qh(\phi(x,t_{*}(x)),a)\}.
\end{eqnarray*}
Consequently, we have that
\begin{align}
w(\phi(x,t_{*}(x))) & = r(\phi(x,t_{*}(x)),\widebar{u}(\phi(x,t_{*}(x))))+Qh(\phi(x,t_{*}(x)),\widebar{u}(\phi(x,t_{*}(x)))).
\label{infmin1}
\end{align}
Combining equations (\ref{infmin1bis}) and (\ref{infmin1}), it follows that 
\begin{eqnarray*}
w(x) & = & -\rho \mathcal{L}_{\alpha}(x,\widehat{u}_{\phi}(w,h)(x)) +L_{\alpha}f(x,\widehat{u}_{\phi}(w,h)(x))+H_{\alpha}r(x,\widehat{u}_{\phi}(w,h)(x)) \nonumber \\
& &+G_{\alpha}h(x,\widehat{u}_{\phi}(w,h)(x)).
\end{eqnarray*}

\bigskip

\noindent
Consider now the case in which $t_{*}(x)=\infty$.
From assumption \ref{Hyp8a}, we obtain that the limit when $t$ tends to infinity of
\begin{align*}
\int_0^t e^{-(\alpha s+\overline{\Lambda}(x,s))} \Bigl[ f(\phi(x,s),\widebar{u}(\phi(x,s))) +\lambda(\phi(x,s),\widebar{u}(\phi(x,s))) Qh(\phi(x,s),\widebar{u}(\phi(x,s))) - \rho \Bigr] ds.
\end{align*}
exists in $\RR\union \{+\infty\}$ and that $w(\phi(x,t))\geq -\rho K_{\xi} - K_h$ for all $t\in[0,+\infty)$.
Therefore, by using equation (\ref{espoir}) we obtain that
\begin{align*}
w(x) \geq -e^{-(\alpha t+\overline{\Lambda}(x,t))} & [\rho K_{\xi} + K_h] +
\int_0^t e^{-(\alpha s+\overline{\Lambda}(x,s))} \Bigl[ f(\phi(x,s),\widebar{u}(\phi(x,s))) \nonumber \\
& +\lambda(\phi(x,s),\widebar{u}(\phi(x,s))) Qh(\phi(x,s),\widebar{u}(\phi(x,s))) - \rho \Bigr] ds,
\end{align*}
and so, the feedback measurable selector $\widehat{u}_{\phi}(w,h) \in \mathcal{S}_{\mathbb{V}}$ satisfies:
\begin{align*}
w(x) \geq & -\rho \mathcal{L}_{\alpha}(x,\widehat{u}_{\phi}(w,h)(x)) + L_{\alpha}f(x,\widehat{u}_{\phi}(w,h)(x)) \nonumber \\
& +\int_0^{t_{*}(x)} e^{-(\alpha s+\overline{\Lambda}(x,s))} \lambda(\phi(x,s),\widebar{u}(\phi(x,s))) Qh(\phi(x,s),\widebar{u}(\phi(x,s)))ds \\
= & -\rho \mathcal{L}_{\alpha}(x,\widehat{u}_{\phi}(w,h)(x)) +L_{\alpha}f(x,\widehat{u}_{\phi}(w,h)(x))+H_{\alpha}r(x,\widehat{u}_{\phi}(w,h)(x)) \\
 & +G_{\alpha}h(x,\widehat{u}_{\phi}(w,h)(x)).
\end{align*}
This shows that
\begin{eqnarray*}
w(x) & = & -\rho \mathcal{L}_{\alpha}(x,\widehat{u}_{\phi}(w,h)(x)) +L_{\alpha}f(x,\widehat{u}_{\phi}(w,h)(x))+H_{\alpha}r(x,\widehat{u}_{\phi}(w,h)(x)) \nonumber \\
& & +G_{\alpha}h(x,\widehat{u}_{\phi}(w,h)(x)).
\end{eqnarray*}

\bigskip

In conclusion, since $\mathbb{V}(x) \subset \mathbb{V}^{r}(x)$ it follows that $\mathcal{R}_{\alpha}(\rho,h)(x)\leq
\mathcal{T}_{\alpha}(\rho,h)(x)$. However, we have shown that $\widehat{u}_{\phi}(w,h) \in \mathcal{S}_{\mathbb{V}}$ satisfies
\begin{eqnarray*}
\mathcal{R}_{\alpha}(\rho,h)(x)& = & -\rho \mathcal{L}_{\alpha}(x,\widehat{u}_{\phi}(w,h)(x)) +
L_{\alpha}f(x,\widehat{u}_{\phi}(w,h)(x))+H_{\alpha}r(x,\widehat{u}_{\phi}(w,h)(x))\nonumber \\
& & +G_{\alpha}h(x,\widehat{u}_{\phi}(w,h)(x)),
\end{eqnarray*}
which is the desired result.
\hfill $\Box$

\bibliography{LongRun}

\end{document}